\def\TheMagstep{\magstep1}      
\def\PaperSize{letter}          
\def\date{25 April 1988}

\let\:=\colon \let\ox=\otimes \let\To=\longrightarrow
\let\vf=\varphi \let\x=\times \let\o=\circ \let\?=\overline
\let\into=\hookrightarrow \let\To=\longrightarrow
\def\onto{\to\mathrel{\mkern-15mu}\to}
 \let\wt=\widetilde
\let\Sum=\sum \def\sum{{\textstyle\Sum}}
\def\IC{{\bf C}} \def\IP{{\bf P}} \def\mY{{\bf m}_Y}
\def\IG{{\bf G}} \def\mG{{\bf m}_{\bf G}} \def\IV{{\bf V}}
\def\pd #1#2{\partial#1/\partial#2}

\def\sdep #1{#1^\dagger}	
\def\and{\hbox{ and }}
\def\Wf{\hbox{\rm W$_f$}}	\def\Af{\hbox{\rm A$_f$}}
\def\IHom{\mathop{{\cal H}\!om}\nolimits}
\def\DONE{*!*}
\def\NextDef #1 {\def\NextOne{#1}%
 \ifx\NextOne\DONE\let\next\relax
 \else\expandafter\xdef\csname#1\endcsname{\TheOp}
  \let\next\NextDef
 \fi \next}
\def\TheOp{\mathop{\rm\NextOne}}
 \NextDef
  Projan Supp Proj Sym Specan Hom cod Ker dist
 *!*
\def\TheOp{{\cal\NextOne}}
\NextDef
  E F G H I J M N O R S
 *!*
\def\TheOp{\hbox{\rm\NextOne}}
\NextDef
 A ICIS
 *!*

\def\item#1 {\par\indent\indent\indent\indent \hangindent4\parindent
 \llap{\rm (#1)\enspace}\ignorespaces}
 \def\inpart#1 {{\rm (#1)\enspace}\ignorespaces}
 \def\part {\par\inpart}

\let\@=@
\catcode`\@=11		

\def\vfootnote#1{\insert\footins\bgroup
 \eightpoint 
 \interlinepenalty\interfootnotelinepenalty
  \splittopskip\ht\strutbox 
  \splitmaxdepth\dp\strutbox \floatingpenalty\@MM
  \leftskip\z@skip \rightskip\z@skip \spaceskip\z@skip \xspaceskip\z@skip
  \textindent{#1}\footstrut\futurelet\next\fo@t}

\def\p.{p.\penalty\@M \thinspace}
\def\pp.{pp.\penalty\@M \thinspace}
\def\(#1){{\rm(#1)}}\let\leftp=(
\def\activeleftp{\catcode`\(=\active}
{\activeleftp\gdef({\ifmmode\let\next=\leftp \else\let\next=\(\fi\next}}

\def\sct#1\par
  {\removelastskip\vskip0pt plus2\normalbaselineskip \penalty-250
  \vskip0pt plus-2\normalbaselineskip \bigskip
  \centerline{\smc #1}\medskip}

\newcount\sctno \sctno=0
\def\sctn{\advance\sctno by 1
 \sct\number\sctno.\quad\ignorespaces}

\def\dno#1${\eqno\hbox{\rm(\number\sctno.#1)}$}
\def\Cs#1){\unskip~{\rm(\number\sctno.#1)}}

\def\proclaim#1 #2 {\medbreak
  {\bf#1 (\number\sctno.#2)}\enspace\bgroup\activeleftp
\it}
\def\endproclaim{\par\egroup\medskip}
\def\pf{\endproclaim{\bf Proof.}\enspace}
\def\lem{\proclaim Lemma } \def\prp{\proclaim Proposition }
\def\cor{\proclaim Corollary }	\def\thm{\proclaim Theorem }
\def\rmk#1 {\medbreak {\bf Remark (\number\sctno.#1)}\enspace}
\def\eg#1 {\medbreak {\bf Example (\number\sctno.#1)}\enspace}

\parskip=0pt plus 1.75pt \parindent10pt
\hsize29pc
\vsize44pc
\abovedisplayskip6pt plus6pt minus2pt
\belowdisplayskip6pt plus6pt minus3pt

\def\TRUE{TRUE}	
\ifx\DoublepageOutput\TRUE \def\TheMagstep{\magstep0} \fi
\mag=\TheMagstep

\newskip\vadjustskip \vadjustskip=0.5\normalbaselineskip
\def\centertext
 {\hoffset=\pgwidth \advance\hoffset-\hsize
  \advance\hoffset-2truein \divide\hoffset by 2\relax
  \voffset=\pgheight \advance\voffset-\vsize
  \advance\voffset-2truein \divide\voffset by 2\relax
  \advance\voffset\vadjustskip
 }
\newdimen\pgwidth\newdimen\pgheight
\def\letter{letter}\def\AFour{AFour}
\ifx\PaperSize\letter
 \pgwidth=8.5truein \pgheight=11truein
 \message{- Got a paper size of letter.  }\centertext
\fi
\ifx\PaperSize\AFour
 \pgwidth=210truemm \pgheight=297truemm
 \message{- Got a paper size of AFour.  }\centertext
\fi

\def\today{\ifcase\month\or	
 January\or February\or March\or April\or May\or June\or
 July\or August\or September\or October\or November\or December\fi
 \space\number\day, \number\year}
\nopagenumbers
 \newcount\pagenumber \pagenumber=1
 \def\advancepagenumber{\global\advance\pagenumber by 1}
\def\folio{\number\pagenum} 
\headline={%
  \ifnum\pagenum=1\firstheadline
  \else
    \ifodd\pagenum\oddheadline
    \else\evenheadline\fi
  \fi
}
\expandafter\ifx\csname date\endcsname\relax \let\dato=\today
	    \else\let\dato=\date\fi
\let\firstheadline\hfill
\def\oddheadline{\eightpoint \rlap{\dato}
 \hfil\headtitle\hfil\llap{\folio}}
\def\evenheadline{\eightpoint\rlap{\folio}
 \hfil\author\hfil\llap{\dato}}
\def\headtitle{\title}

 \newdimen\fullhsize \newbox\leftcolumn
 \def\fulline{\hbox to \fullhsize}
\def\doublepageoutput
{\let\lr=L
 \output={\if L\lr
           \global\setbox\leftcolumn=\columnbox \global\let\lr=R%
          \else \doubleformat \global\let\lr=L
	  \fi
        \ifnum\outputpenalty>-20000 \else\dosupereject\fi
	}%
 \def\doubleformat{\shipout\vbox{%
     \ifx\PaperSize\AFour
	   \fulline{\hfil\box\leftcolumn\hfil\columnbox\hfil}%
     \else
	   \fulline{\hfil\hfil\box\leftcolumn\hfil\columnbox\hfil\hfil}%
     \fi             }%
     \advancepageno
}
 \def\columnbox{\vbox
   {\ifvoid255\headline={\hfil}\nopagenumbers\else
    \makeheadline\pagebody\makefootline\advancepagenumber\fi}%
   }%
\fullhsize=\pgheight \hoffset=-1truein
 \voffset=\pgwidth \advance\voffset-\vsize
  \advance\voffset-2truein \divide\voffset by 2
  \advance\voffset\vadjustskip
 
}
\ifx\DoublepageOutput\TRUE \let\pagenum=\pagenumber\doublepageoutput
 \else \let\pagenum=\pageno \fi

 \font\twelvebf=cmbx12		
 \font\smc=cmcsc10		

\def\eightpoint{\eightpointfonts
 \setbox\strutbox\hbox{\vrule height7\p@ depth2\p@ width\z@}%
 \eightpointparameters\eightpointfamilies
 \normalbaselines\rm
 }
\def\eightpointparameters{%
 \normalbaselineskip9\p@
 \abovedisplayskip9\p@ plus2.4\p@ minus6.2\p@
 \belowdisplayskip9\p@ plus2.4\p@ minus6.2\p@
 \abovedisplayshortskip\z@ plus2.4\p@
 \belowdisplayshortskip5.6\p@ plus2.4\p@ minus3.2\p@
 }
\newfam\smcfam
\def\eightpointfonts{%
 \font\eightrm=cmr8 \font\sixrm=cmr6
 \font\eightbf=cmbx8 \font\sixbf=cmbx6
 \font\eightit=cmti8
 \font\eightsmc=cmcsc8
 \font\eighti=cmmi8 \font\sixi=cmmi6
 \font\eightsy=cmsy8 \font\sixsy=cmsy6
 \font\eightsl=cmsl8 \font\eighttt=cmtt8}
\def\eightpointfamilies{%
 \textfont\z@\eightrm \scriptfont\z@\sixrm  \scriptscriptfont\z@\fiverm
 \textfont\@ne\eighti \scriptfont\@ne\sixi  \scriptscriptfont\@ne\fivei
 \textfont\tw@\eightsy \scriptfont\tw@\sixsy \scriptscriptfont\tw@\fivesy
 \textfont\thr@@\tenex \scriptfont\thr@@\tenex\scriptscriptfont\thr@@\tenex
 \textfont\itfam\eightit	\def\it{\fam\itfam\eightit}%
 \textfont\slfam\eightsl	\def\sl{\fam\slfam\eightsl}%
 \textfont\ttfam\eighttt	\def\tt{\fam\ttfam\eighttt}%
 \textfont\smcfam\eightsmc	\def\smc{\fam\smcfam\eightsmc}%
 \textfont\bffam\eightbf \scriptfont\bffam\sixbf
   \scriptscriptfont\bffam\fivebf	\def\bf{\fam\bffam\eightbf}%
 \def\rm{\fam0\eightrm}%
 }

 \newcount\refno \refno=0	 \def\NoKey{*!*}
 \def\MakeKey{\advance\refno by 1 \expandafter\xdef
  \csname\TheKey\endcsname{{\number\refno}}\NextKey}
 \def\NextKey#1 {\def\TheKey{#1}\ifx\TheKey\NoKey\let\next\relax
  \else\let\next\MakeKey \fi \next}
 \def\RefKeys #1\endRefKeys{\expandafter\NextKey #1 *!* }
\def\SetRef#1 #2,#3\par{%
 \hang\llap{[\csname#1\endcsname]\enspace}%
  \ignorespaces{\smc #2,}
  \ignorespaces#3\unskip.\endgraf
 }
 \newbox\keybox \setbox\keybox=\hbox{[8]\enspace}
 \newdimen\keyindent \keyindent=\wd\keybox
\def\references{
  \bgroup   \frenchspacing   \eightpoint
   \parindent=\keyindent  \parskip=\smallskipamount
   \everypar={\SetRef}}
\def\endreferences{\egroup}

 \def\serial#1#2{\expandafter\def\csname#1\endcsname ##1 ##2 ##3
  {\unskip\ #2 {\bf##1} (##2), ##3}}
 \serial{ajm}{Amer. J. Math.}
 \serial{ann}{Annals Math.}
 \serial{asens}{Ann. Scient. \'Ec. Norm. Sup.}
 \serial{cmh}{Comment. Math. Helvetici}
 \serial{comp}{Compositio Math.}
 \serial{conm}{Contemp. Math.}
 \serial{crasp}{C. R. Acad. Sci. Paris}
 \serial{dlnpam}{Dekker Lecture Notes in Pure and Applied Math.}
 \serial{faa}{Funct. Anal. Appl.}
 \serial{invent}{Invent. Math.}
 \serial{ja}{J. Algebra}
 \serial{jag}{J. Alg. Geom.}
 \serial{ma}{Math. Ann.}
 \serial{mpcps}{Math. Proc. Camb. Phil. Soc.}
 \serial{splm}{Springer Lecture Notes in Math.}
 \serial{tams}{Trans. Amer. Math. Soc.}

\def\UThin{\penalty\@M \thinspace\ignorespaces}
\def\relaxnext@{\let\next\relax}
\def\cite#1{\relaxnext@
 \def\nextiii@##1,##2\end@{\unskip\space{\rm[\SetKey{##1},\let~=\UThin##2]}}%
 \in@,{#1}\ifin@\def\next{\nextiii@#1\end@}\else
 \def\next{{\rm[\SetKey{#1}]}}\fi\next}
\newif\ifin@
\def\in@#1#2{\def\in@@##1#1##2##3\in@@
 {\ifx\in@##2\in@false\else\in@true\fi}%
 \in@@#2#1\in@\in@@}
\def\SetKey#1{{\bf\csname#1\endcsname}}

\catcode`\@=12 

\RefKeys
  BMM BPW B-R Fulton G-1 G-2 G-3 G-4 GaM Gr G-M G H&M82 H&M83 H-M H-M-S
Hir R-K Kl98 K-T KT97 L L-R L-S LT81 LT88 LJT Lip Lo84 N80 P R87 R89 T-1
T75 T76 T77 T1.5 T81 T-2 Trot Verdier
 \endRefKeys

\def\author{T. Gaffney and S. Kleiman}
\def\title{Specialization of integral dependence for modules}

\leavevmode
 \bigskip\bigskip
 \centerline{\twelvebf \title}
 \bigskip
 \footnote{}{\noindent 1991 {\it Mathematics Subject Classification}.
   Primary 32S15; Secondary 14B05, 13H15.}
 \footnote{}{It is a pleasure to thank Robert Gassler for reading over
the manuscript carefully and making a number of helpful comments.}%
 \centerline{\smc Terence GAFFNEY\footnote{$^{1}$}{%
    Supported in part by NSF grant 9403708-DMS.}
 and Steven L. KLEIMAN\footnote{$^{2}$}{%
    Supported in part by NSF grants 9106444-DMS and 9400918-DMS.}
            }
\vskip15pt plus12pt minus12pt
{\parindent=24pt \narrower \noindent \eightpoint
 {\smc Abstract.}\enspace
 We establish the principle of specialization of integral dependence for
submodules of finite colength of free modules, as part of the general
algebraic-geometric theory of the Buchsbaum--Rim multiplicity.  Then we
apply the principle to the study of equisingularity of ICIS germs,
obtaining results for such equisingularity conditions as Whitney's
Condition A, Thom's Condition \Af, and the Relative Whitney Condition
\Wf.  Notably, we describe these conditions for analytic families in
terms of various numerical invariants, which, for the most part, depend
only on the members of a family, not on its total space.
 \par }

\sct Introduction

Describing the structure of a singular set remains a basic, but elusive,
goal of complex-analytic geometry.  If the set is a member of an
analytic family, then it is often easier to tell when the set's
structure is similar to that of the general member.  In this paper, we
study analytic families of germs of isolated complete-intersection
singularities, or {\bf ICIS} germs, as a step toward the general study.
We develop some algebraic tools and a geometric point of view that
enable us to describe many equisingularity conditions for these families
in terms of numerical invariants.  The invariants, for most conditions,
depend {\it only on the individual members\/} of a family, not on its
total space.

The basic numerical invariants we use are certain Buchsbaum--Rim
multiplicities.  They arise from the column space of the Jacobian matrix
of a given ICIS germ.  This column space is known as the {\it Jacobian
module}.  It is, in a natural way, a submodule of finite colength in a
free module over the local ring of the germ, and the associated
invariants govern its integral closure.  In this context, our main
algebraic theorem is a generalization, from ideals to modules, of
Teissier's ``principle of specialization of integral dependence,''
\cite{T-1, 3.2, p.~330} and \cite {T1.5, App.~I}.  We prove the theorem
and some related results in Sections~1 and 2.  In Section~3, we treat
the related notion of strict dependence, which we use to handle those
equisingularity conditions that require a little more than integral
dependence to ensure that they hold.

In Sections~4, 5, and 6, we apply the results of Sections~1, 2 and 3 to
the study of various equisingularity conditions on families of ICIS
germs.  More specifically, in Section~4, we study Whitney's Condition~A.
We prove notably that A is satisfied at the origin if every associated
multiplicity of the Jacobian module is constant across the family.

In Section~5, we fix a function $f$ on the total space $X$, and study
{\it Thom's Condition\/} \Af.  Notably, we generalize a celebrated
theorem of L\^e and Saito \cite{L-S}: using Parame\-swar\-an's
construction in \cite{P}, we relate \Af\ to information about vanishing
cycles as encoded in a sequence of Milnor numbers.  In addition, we
refine a theorem of Brian\c con, Maisonobe and Merle's \cite{BMM,
Thm.~4.2.1, p.~541} in the present setting of a family of ICIS germs:
their theorem requires the fulfillment of Condition A and the stratified
local topological triviality of every linear projection to the singular
locus, whereas ours requires only the constancy of a Buchsbaum--Rim
multiplicity, or equivalently, of two Milnor numbers.

Finally, in Section 6, we study \Wf, the standard relative form of
Whitney's Condition~B.  Notably, we establish three necessary and
sufficient conditions for \Wf\ to hold.  Two are {\it memberwise\/}
conditions: the constancy of two sequences of Milnor numbers, and the
constancy of a single Buchsbaum--Rim multiplicity.

The remaining condition is global: denote the locus of central points by
$Y$, and set $Z:=f^{-1}0$; then both pairs, $(X-Y,Y)$ and $(Z-Y,Y)$ must
satisfy Whitney's Condition B at the origin.  The sufficiency of this
condition was established by Brian\c con, Maisonobe and Merle in
\cite{BMM, Thm.~4.3.2, p.~543} in a more general setting, and it is
recovered in the present setting via a new proof.  The new proof
illustrates the use of integral closure methods, and lays the foundation
for further progress in the study of \Wf\ both for nonisolated
singularities and for families of isolated singularities.

Some of the present work, together with its extension in \cite{GaM}, is
explained and developed by the second author in \cite{Kl98}.  See also
\cite{KT97, (1.7)}.

Let us now discuss the contents in more detail, stressing the philosophy
of our approach.  Let $(X,0)\to(Y,0)$ be a map of germs of complex
analytic spaces.  (As is conventional when dealing with germs, we often
let it go without saying that any analytic set, even $\IC^n$ itself,
should be {\it replaced\/} by a suitably small neighborhood of the
central point, or ``origin,'' whenever appropriate.)  Let $X(y)$ denote
the fiber over the point $y\in Y$.  Assume that the $X(y)$ are
equidimensional of the same dimension $d$, where $d\ge1$, that $X$ is
equidimensional, and that $Y$ has dimension at least 1.  Let
$\E:=\O_X^p$ be a free module of rank $p$ at least 1.  Let $\M$ be a
coherent submodule such that, on each fiber $X(y)$, the image of the
restriction of $\M$ in that of $\E$ has finite colength.

Theorem~(1.8) is our generalization of Teissier's principle.  Our
theorem gives a criterion, valid for any rank $p$, for a section $h$ of
$\E$ to be integrally dependent on $\M$, or equivalently, for the
submodule $\H$ generated by $h$ and $\M$ to lie in the integral closure
of $\M$ (in other words, for $\M$ to be a {\it reduction\/} of $\H$).
As we show in Sections~5 and 6, Theorem~(1.8) is a powerful tool for
establishing sufficient conditions for an equisingularity condition to
hold.

Since Theorem~(1.8) is central to our work, let us outline its proof.
Let $\S\E$ denote the symmetric algebra, and $\R\H$ and $\R\M$ the Rees
algebras, which are the subalgebras generated by $\H$ and $\M$.  Form
the analytic homogeneous spectra,
  $$P:=\Projan(\S\E),\ Q:=\Projan(\R\H)\and P':=\Projan(\R\M).$$
 Then the inclusion of $\M$ into $\H$ induces a well-defined finite map
from $Q$ to $P'$ if and only if $\M$ is a reduction of $\H$;
cf.~\cite{K-T, (2.6)}.

 In Sections 8 and 9 of \cite{K-T}, this setup is used to associate
various multiplicities to the pair $(\M,\H)$.  It is shown there, in
Theorem~(9.5) and Corollary~(9.7), that, if $\M$ is {\it not\/} a
reduction of $\H$ everywhere, but is so away from the origin, then
certain of these multiplicities are positive.  In turn, by
Corollary~(10.2) of~\cite{K-T}, this positivity implies, surprisingly
enough, that the dimension of the central fiber of $P'/X$ is of maximal
dimension, namely, $n-1$.  (The recent paper \cite{KT97} gives, in
(1.4), a substantially shorter, simpler, and more direct proof that this
dimension is maximal.)  Thus, if we can control $P'$, then we can force
$\H$ to be in the integral closure of $\M$.

For each $y$, let $e(y)$ denote the {\it Buchsbaum--Rim multiplicity} of
the pair induced on $X(y)$ by $(\E,\M)$; so $e(y)/(d+p-1)!$ is the
coefficient of $n^{d+p-1}$ in the polynomial whose value for $n\gg0$ is
the (vector space) dimension of $(\S_n\E/\R_n\M)(y)$.  Just as Teissier
did in the case of an ideal, we show that the constancy of $e(y)$
implies the existence of a reduction of $\M$ with the minimal number of
generators, and therefore yields the desired upper bound on the
dimension of the central fiber of $P'/X$.  The proof depends heavily on
the upper semicontinuity of $e(y)$, which we establish in
Proposition~(1.1).  Thus Theorem~(1.8) is proved.

We can also study the integral dependence of modules by reducing this
study to that of the integral dependence of ideals; compare with
\cite{R89} and \cite{H-M, p.~160}.  Indeed, in the setup above, let
$\rho(\M)$ be the ideal on $P$ associated to the homogeneous ideal in
the symmetric algebra $\S\E$ generated by $\M$ viewed in degree 1.  By
Proposition~(3.4), the module $\H$ is integrally dependent on the module
$\M$ if and only if the ideal $\rho(\H)$ is integrally dependent on the
ideal $\rho(\M)$.

An important difference between the case of modules and that of ideals
shows up when we form the blowup $B$ of $P$ with respect to $\rho(\M)$.
Namely, if an ideal $I$ has finite colength, then the corresponding
exceptional divisor is naturally a projective scheme, and its degree is
simply the multiplicity of $I$.  On the other hand, even though $\M$ has
finite colength in $\E$, its associated ideal $\rho(\M)$ may define a
subset of $P$ of positive dimension, and then the exceptional divisor of
$B$ is naturally a biprojective scheme.  So it has a series of
bidegrees, its {\it Segre numbers}, defined by intersecting it with the
various powers of the two hyperplane classes.  The partial sums of these
Segre numbers are the {\it associated multiplicities\/} of $\M$.  They
can also be computed directly from a set of generators of $\M$.

In Section 2, we advance the theory of these Segre numbers and
associated multiplicities, relating them to the dimension of the
exceptional divisor of $B$.  By using these multiplicities instead of
the multiplicity of an ideal, we can generalize to ICIS germs many
results about families of hypersurface germs.  We do so in Section~4 for
some results about Whitney's Condition A.

Section~3 concerns strict dependence.  In the case of an ideal $\I$,
Lejeune and Teissier \cite{LJT, pp.~46--48} showed that, if $\I$ has
finite colength, then a germ $h$ is strictly dependent on $\I$ if and
only if, on every component of the exceptional divisor in the normalized
blowup of $\I$, the order of vanishing of the pullback of $h$ is
strictly greater than the order of vanishing of the pullback of $\I$.
Proposition (3.5) gives the generalization of this result to the case of
a module $\M$; instead of blowing up $\I$, we blow up $\rho(\M)$.

In the case where $\I$ defines a family of ideals of finite colength and
constant multiplicity, Teissier\cite {T1.5, App.~I} used the constancy
to force the components of the exceptional divisor to surject onto the
parameter space.  Then it is easy to show that, if the restriction of a
germ $h$ to a general fiber is strictly dependent on the restriction of
$\I$, then the restriction of $h$ to every fiber is strictly dependent.
In the module case, using the associated multiplicities, we make a
similar argument and draw a similar conclusion in the course of proving
Proposition~(4.2).

In Sections~4, 5, and 6, we assume that the fibers $X(y)$ represent ICIS
germs, and abusing notation, we also let $Y$ denote the locus of central
points.  We embed $(X,0)$ in $(\IC^n,0)$ so that $(Y,0)$ is the germ of
a linear subspace.  In Sections~5, and 6, we also consider a function
germ $f\:X\to \IC$ and its zero set $Z$, and we assume $Z\supset Y$.

Section 4 concerns {\it Whitney's Condition\/} A.  Let $S$ denote the
singular locus of $X$, and assume that $S$ is smooth.  By definition,
Condition~A is satisfied by the pair $(X-S,S)$ at $0\in S$ if the
tangent space $T_0S$ lies in every hyperplane obtained as a limit of
hyperplanes, each tangent to $X-S$ as the point of contact approaches
$0$.

Theorem~(4.2) asserts that Condition A is satisfied if every associated
multiplicity of the Jacobian module of the fiber $X(y)$ is constant in
$y$.  The theorem is illustrated in Examples~(4.3) and (4.4).  In the
first, we work out, from our viewpoint, Trotman's example, showing that
no fiberwise criterion for Condition A can be necessary as well as
sufficient.  In the second, we consider a case where $S$ is larger than
$Y$.

Section~5 concerns {\it Thom's Condition\/} \Af.  Let $\Sigma(f)$ denote
the union of the singular points of the fibers of $f$.  By definition,
$(X-\Sigma(f),Y)$ satisfies \Af\ at $y\in Y$ if the tangent space to $Y$
at $y$ lies in every limit tangent hyperplane at $y$ to the fibers of
$f$.  In Lemma~(5.1), we connect \Af\ to the theory of integral closure
via the augmented Jacobian module.

Thom introduced \Af\ as the primary condition guaranteeing the local
topological triviality of the family of functions defined by $f$.
Condition \Af\ is also important because of its relationship (which is
well understood in only a few cases) to the vanishing cycles.  For
example, in a recent paper \cite{G-M}, Green and Massey show that
certain information about the vanishing cycles implies \Af\ for families
with generalized isolated singularities.  In the case where $X=\IC^n$
and $\Sigma(f)=Y$, L\^e and Saito \cite{L-S} showed that \Af\ is implied
by the constancy of the Milnor number.

Theorem~(5.2) generalizes the L\^e--Saito theorem to the case where $X$
is a complete intersection.  This generalization is based on the
following construction of Parameswaran's in \cite{P}, which reduces the
study of a family of ICIS germs to the study of an isolated singularity
defined by the vanishing of a single function inside an ambient space
with an isolated singularity.  Given an ICIS germ $(X,0)$ with embedding
codimension $k$, for each $i$ with $0\le i\le k$ let $\mu _i$ be the
smallest of all the Milnor numbers of ICIS germs that serve as total
spaces of $i$-parameter (flat) deformations of $(X,0)$.  Denote this
sequence of numbers by $\mu_*$.

A chain of ICIS germs, each of codimension $1$ inside the next, whose
sequence of Milnor numbers is $\mu _*$ is said to be $\mu _*$-{\it
minimal}.  Given the germ of a family $X/Y$ with $Y$ smooth,
Parameswaran constructed a chain of deformations of the family such that
the parameter spaces are smooth over $Y$ and the central chain is $\mu
_*$-minimal; moreover, if the $\mu _*$-sequence is constant in the given
family, then the chain of deformations is $\mu _*$-minimal.

Given a chain of deformations of $(X,0)$ that is $\mu _*$-minimal, we
can look, on the deformation $X_i$, at the function $f_{i-1}$ that
defines the deformation $X_{i-1}$, and ask that A$_{f_{i-1}}$ be
satisfied at 0 by the pair $(X_i-Y,Y)$ where $Y$ is the common singular
locus of all the $X_i$.  Theorem~(5.2) gives a necessary and sufficient
condition for A$_{f_{i-1}}$ to hold for every $i$, namely, the constancy
of $\mu_*$.

Indeed, this constancy turns out to be equivalent to the constancy of
the multiplicity of the relative augmented Jacobian module of $X_i$ and
$f_{i-1}$; the principle of specialization of integral dependence then
shows that the generators of the relative Jacobian module that come from
the partial derivatives along $Y$ are dependent on the augmented module.
The original theorem of L\^e and Saito, which is reproved, shows that
this dependence is strict at the top of the chain, where the ambient
space is just the affine space.  It remains to push the strict
dependence down the chain.  A careful examination of the proof shows a
similarity in the role played by the elements in the chain and the
associated multiplicities.

We prove the sufficiency as follows.  Because the number of generators
of the relative augmented Jacobian submodule is always the minimum
number needed to generate a module of finite nonzero colength, it
follows from Proposition~(1.5)(3) that, if A$_{f_{i-1}}$ holds, then the
multiplicity of the augmented Jacobian module of the fibers of $X_i$ and
$f_{i-1}$ is independent of the parameter value; hence, so is the
$\mu_*$-sequence.

We also refine, in our case of a family of ICIS germs, a theorem proved
by Brian\c con, Maisonobe and Merle \cite{BMM, Thm.~4.2.1, p.~541}.
Recalled before our Theorem~(5.3), their theorem essentially asserts this:
\Af\ holds along $Y$ if, for every linear retraction $r$ to $Y$, the
restriction $r|(X,Z,Y)$ is stratified locally topologically trivial, and
if both pairs $(X-Y,Y)$ and $(Z-Y,Y)$ satisfy Whitney's Condition A along
$Y$.  Our machinery allows us to prove, in Theorem~(5.3), that \Af\
holds at $0$ assuming only, for every $r$, that the Buchsbaum--Rim
multiplicity of the augmented Jacobian module of the fibers of $X$ and
$f$ is defined and is independent of the parameter value, or
equivalently that, the germs of the fibers of the restrictions $r|X$ and
$r|Z$ have isolated singularities, and their Milnor numbers are
independent of the parameter value.

Our work leads us to conjecture that \Af\ holds whenever the
Buchsbaum--Rim multiplicity is defined and is independent of the
parameter value, or the Milnor numbers are so.  After the present work
was completed, Massey and the first author \cite{GaM, (5.8)} proved this
conjecture for families of ICIS germs via a careful study of the
conormal variety; see also \cite{Kl98} and \cite{KT97, (1.7)}.  At the
moment, this multiplicity is defined only when $X$ is the total space of
a family of ICIS germs.  However, we conjecture that, once the theory of
multiplicity has been extended to cover modules of infinite colength,
then the independence of the multiplicity of the augmented Jacobian
module will always be equivalent to \Af.

Assume also that $(Z-Y,Y)$ satisfies Whitney's Condition A at the
origin.  Then \Af\ holds if the Milnor numbers of the fibers of the two
restrictions $r|X$ and $r|Z$ are independent of the parameter value for
only a single $r$; see Theorem~(5.5).  We prove it via a close analysis
of the relative conormal space of $f$ using the principle of
specialization of integral dependence to gain control.  The hypothesis
on $Z$ gives additional information on the relative conormal space,
since it contains the conormal space of $Z$.

 Section 6 treats the condition \Wf.  By definition, it is satisfied by
the pair $(X-\Sigma(f),Y)$ at $y\in Y$ if each tangent plane to the
fiber of $f$ at an $x$ in $X-\Sigma(f)$ approaches the tangent plane of
$Y$ at $y$ as fast as $x$ approaches $Y$.  In Proposition~(6.1), we
connect \Wf\ to the theory of integral closure via the augmented
Jacobian module again; we also recover the integral closure condition of
L\^e and Teissier \cite{LT88, Prop.~1.3.8} between ideals on the
relative conormal space.  In Lemma~(6.2) we illustrate the usefulness of
our integral closure condition by using it to give a new proof of a
basic transversality result of Henry and Merle.

It is natural to ask if there is a numerical invariant of the fibers of
$X$ and $f$ over $Y$ whose constancy ensures that \Wf\ is satisfied at
0.  An affirmative answer is given in Theorem~(6.4): a suitable
invariant is the multiplicity $em(y)$ of the product of the maximal
ideal of $X(y)$ and its augmented Jacobian ideal.  In fact, assuming
that each fiber $Z(y)$ has an isolated singularity at $0$, we prove that
the constancy of $em(y)$ is both necessary and sufficient for \Wf\ to be
satisfied.  The theorem also proves another necessary and sufficient
condition: that $(X-Y,Y)$ and $(Z-Y,Y)$ satisfy both Whitney conditions
at 0.  The key ingredient in our proof is again Theorem~(1.8).

We interpret $em(y)$ topologically in Lemma~(6.3); namely, $em(y)$ is
equal to a linear combination with certain binomial coefficients of the
sum of the Milnor numbers of the plane sections of $X(y)$ and $Z(y)$.
Hence the constancy of $em(y)$, and so that of \Wf, is equivalent to the
constancy of the Milnor numbers of these plane sections.  This
equivalence is also part of Theorem~(6.4).

\sctn Specialization of integral dependence

 Let $F\:(X,x_0)\to(Y,y_0)$ be a map of germs of complex analytic
spaces, which need not be reduced.  Assume that the fibers $X(y)$ are
equidimensional of the same dimension $d$ at least 1 and that $Y$ has
dimension at least 1.  Let $\E:=\O_X^p$ be a free module of rank $p$ at
least 1, and set $r:=d+p-1$.  Let $\M$ be a coherent submodule of $\E$.
Set $S:=\Supp(\E/\M)$, and assume that $S$ is finite over $Y$.  Finally,
for each $y\in Y$, denote the {\it Buchsbaum--Rim multiplicity} of the
pair that $(\E,\M)$ induces on $X(y)$ by $e(y)$.

Let $\N$ be a coherent submodule of $\M$.  Let $\S\E$ denote the
symmetric algebra, and $\R\M$ and $\R\N$ the subalgebras
generated by $\M$ and $\N$; say
        $$\textstyle\S\E=\bigoplus_n\S_n\E,\ \R\M=\bigoplus_n\R_n\M
        \hbox{ and }\R\N=\bigoplus_n\R_n\N$$
 are the decompositions into graded pieces.  Form the analytic
homogeneous spectra,
  $$P:=\Projan(\S\E),\ P':=\Projan(\R\M)\and P'':=\Projan(\R\N).$$
 Recall that, if $\R\M$ is a finitely generated
$\R\N$-module, then $\N$ is called a {\it reduction\/} of $\M$, and the
sections of $\M$ are said to be {\it integrally dependent} on $\N$.  A
different, but equivalent, definition of integral dependence is
discussed at the beginning of Section~3 (and a third definition is
mentioned there in passing).

If $\N$ is a reduction of $\M$, then the following three conditions
obtain:
 \smallskip\item i {\it
 on each fiber $X(y)$ the Buchsbaum--Rim multiplicity arising from
 $(\E,\N)$ is defined and equal to the multiplicity $e(y)$ arising from
 $(\E,\M)$;}
 \item ii {\it
 $\N$ is equal to $\E$ at every point $x$ of $X-S$, that is,
$\Supp(\E/\N)\subset S$;}
 \item iii {\it
  the inclusion $\R\N\into\R\M$ induces a finite surjective map $P'\onto
P''$, which is an isomorphism off $S$.}
 \smallskip\noindent
 Condition (i) obtains by \cite{K-T, (6.7a)(iii)(a), \p.204}.  Condition
(ii) obtains by \cite{K-T, (2.4), \p.182}.  The argument is simple.  By
Nakayama's lemma, it suffices to note that the map of fibers
$\N(x)\to\E(x)$ is surjective.  However, its image is a vector subspace
of $\E(x)$; so, if a basis of the former is extended to a basis of the
latter, then
  the image of $\R\N(x)$ in the polynomial ring $\S\E(x)$ is the
subring generated by a subset of variables; yet the larger ring is a
finitely generated module over the smaller one.  Condition (iii) is
clearly not only necessary, but also sufficient, for $\N$ to be a
reduction of $\M$.

The results below are well known in the case $p=1$; see Lipman's masterful
treatment \cite{Lip}, for example.  For arbitrary $p$, Proposition\Cs5)
relates the existence of a reduction generated by $r$ elements to the
constancy of $e(y)$.  Its proof involves Proposition\Cs1) and
Lemmas\Cs2) and\Cs4).  Proposition\Cs1) asserts the upper
semicontinuity of $e(y)$.  Lemma\Cs2) asserts that if $\N$ is a
reduction of $\M$ fiberwise, then it is so globally over a dense open
subset of $Y$.  Example\Cs3) shows the necessity of passing to an open
subset.  Lemma\Cs4) gives a geometric criterion for $\N$ to be a
reduction of $\M$.  Remark\Cs6) suggests that part of Proposition\Cs5)
should hold in greater generality.  Lemma\Cs7) gives a necessary and
sufficient geometric condition for the existence of a reduction
generated by $r$ elements.  Finally, Theorem\Cs8) rests on all the
preceding results; it gives a generalization of Teissier's principle of
the specialization of integral dependence \cite{T1.5, App.~1}.

\prp 1
 The function $y\mapsto e(y)$ is Zariski upper semicontinuous.
 \pf
 For each $y$, the number $e(y)/r!$ is the coefficient of $n^r$ in the
polynomial whose value is eventually the (vector space) dimension,
        $$\lambda(n,y):=\dim\,\bigl(F_*(\S_n\E/\R_n\M)\bigr)(y).$$
 The polynomial has degree $r$, or else vanishes.  Its value is
$\lambda(n,y)$ for all $n$ at least $n_0$, where $n_0$ depends on $y$.
(See \cite{B-R, bot. \p.213} or \cite{K-T, (5.10)(i)(ii),
\pp.199--200}.)

For each $n$, the $\O_Y$-module $F_*(\S\E_n/\R_n\M)$ is coherent;
hence, $y\mapsto \lambda(n,y)$ is upper semicontinuous.  Therefore,
$y\mapsto e(y)$ is `nondecreasing'; that is, if $A$ is a Zariski
closed irreducible subset of $Y$ and if $\eta$ is its generic point,
then $e(\eta)\le e(y)$ for all $y$ in $A$.  However, it is less obvious
that $A$ contains a nonempty Zariski relatively open subset $U$,
independent of $n$, on which $y\mapsto \lambda(n,y)$ is constant.  To
prove it, we may replace $Y$ by $A$, given its reduced structure, and
replace $X$, $\M$, and so forth by their restrictions.

Form the bigraded $\O_X$-algebra $\R\M\ox\S\E$ and its bigraded module
        $$\F:=\bigoplus\nolimits_{p\ge0,q\ge1} \F_{p,q} \hbox{ where }
        \F_{p,q}:=\R_p\M\S_q\E/\R_{p+1}\M\S_{q-1}\E.$$
 Clearly, $\F$ is generated by $\F_{0,1}$, which is equal to
$\S_1\E/\R_1\M$.  Therefore,  $\F$ is an $\O_S$-module.  So $F_*\F$ is a
finitely generated module over $F_*((\R\M\ox\S\E)|S)$, which is a finitely
generated bigraded $\O_Y$-algebra.  Therefore, by the lemma of generic
flatness, there is a nonempty open subset $U$ of $Y$ on which $F_*\F$
is flat.  Hence, on $U$, each $F_*\F_{p,q}$ is flat.

It follows that, on $U$, the formation of each $F_*\F_{p,q}$ commutes
with restriction to the fibers.  Indeed, since $F|S$ is finite,
$(F_*\F_{p,q})(y)$ is equal to $F_*(\F_{p,q}(y))$ for any $y\in Y$.
Moreover, if $y\in U$, then the formation of $\F_{p,q}$ commutes with
restriction to $X(y)$; this claim will hold, clearly, if the natural
map,
        $$(\R_p\M\S_q\E)(y)\To \S_{p+q}\E(y)$$
  is injective.  It is trivially injective if $p=0$.  Proceeding by
induction on $p$, consider the short exact sequence,
        $$0\To\R_{p+1}\M\S_{q-1}\E\To\R_p\M\S_q\E\To\F_{p,q}\To0.$$
Since $F_*\F_{p,q}$ is flat over $U$ and since $F|S$ is finite, also
$\F_{p,q}$ is flat over $U$. Therefore, the first of these two maps,
        $$(\R_{p+1}\M\S_{q-1}\E)(y)\To(\R_p\M\S_q\E)(y)\To\S_{p+q}\E(y),$$
 is injective.  The second is injective by induction.  So the
composition is injective, as required.

Work on $U$.  Since the formation of $F_*\F_{p,q}$ commutes with
restriction, clearly
        $$\sum_{p+q=n}\dim (F_*\F_{p,q})(y)=\lambda(n,y).$$
 Each $F_*\F_{p,q}$ is flat, so locally free; hence, each function
$y\mapsto\dim F_*\F_{p,q}(y)$ is constant.  Therefore, $y\mapsto
\lambda(n,y)$ is constant.  Thus the proposition is proved.

\lem 2
 Assume that there is a dense Zariski open subset $V$ of $Y$ such
that, for each $y$ in $V$, the image in $\E(y)$ of $\N$ is a reduction
of the image of $\M$.  Then there is a smaller dense Zariski open
subset $U$ of $Y$ over which $\N$ is a reduction of $\M$.
 \pf
 Clearly, it suffices to find a dense Zariski open subset $U$ of $V$
and an integer $k$ such that the inclusion map,
        $$\N\R_k\M\to\R_{k+1}\M,\dno2.1$$
 is surjective over $U$.  By Nakayama's lemma, we may assume that $Y$ is
reduced.  Then there is a dense Zariski open subset $U$ of $Y$ such that
the restriction,
        $(\R_k\M)(y)\to\S_k\E(y)$,
 is injective for all $k$ and for all $y$ in $U$; the existence of $U$
was established in the proof of Proposition\Cs1).  Replace $U$ by $U\cap
V$.

By hypothesis, for each $y$ in $U$, there exists a $k$ such that the
image of the composition,
        $$\N(y)(\R_k\M)(y)\to(\R_{k+1}\M)(y)\to\S_{k+1}\E(y),$$
 is equal to the image of the second map.  Since the second map is
injective, the first map is surjective.  Hence, Nakayama's lemma implies
that, at each point of the fiber $X(y)$, the map \Cs2.1) is surjective.
Therefore, $X(y)$ is contained in the maximal open set on which \Cs2.1)
is surjective, namely, the complement of the support $S_k$ of the
cokernel of \Cs2.1).

On the other hand, $S_k\subset S$; in other words, \Cs2.1) is
surjective at every $x$ off $S:=\Supp(\E/\M)$, as we'll now see.  Set
$y:=F(x)$.  By hypothesis, the image in $\E(y)$ of $\N$ is a reduction
of the image of $\M$.  Since $x\notin S$, the image of $\M$ is equal to
$\M(y)$ at $x$.  Hence, the image of $\N$ is equal to $\M(y)$ at $x$
because of Condition (ii) recalled at the beginning of this section.
In other words, the map $\N\to\M(y)$ is surjective at $x$.  Therefore,
by Nakayama's lemma, the inclusion $\N\to\M$ is surjective at $x$.  So
\Cs2.1) is surjective at $x$, as claimed.

By hypothesis, $S$ is finite over $Y$.  Hence, since $S_k\subset S$,
the image $A_k$ of $S_k$ is a closed analytic subset of $Y$.  Pick a
point $y$ in each component of $U$, and pick a $k$ large enough so that
$A_k$ contains none of these $y$.  Then $U-A_k$ a dense Zariski open
subset $U$ of $V$ on which \Cs2.1) is surjective, as required.

 \eg 3 The open subset provided by Lemma\Cs2) may have to be strictly
smaller than the given open subset (in other words, fiberwise integral
dependence does not imply dependence at the level of the total space).
For example, let $X$ be the $(s,t)$-plane, $Y$ the $s$-line, and
$F\:X\to Y$ the projection.  Let $\E:=\O_X$, let $\N:=(st^2+t^3,
st^4)$, and let $\M:=(\N,t^3)$.  Then, on each fiber of $X/Y$, the
ideals induced by $\M$ and $\N$ are equal.  However, $\N$ is not a
reduction of $\M$ (in other words, $t^3$ is not integrally dependent on
$\N$).  Indeed, otherwise, under the map from the $u$-line into $X$
given by $u\mapsto (u,-u)$, the ideals $\M$ and $\N$ would induce
two ideals, where the second is a reduction of the first; however, $\M$
and  $\N$ induce  $(u^3)$ and  $(u^5)$, and it is easy to see that the
latter is not a reduction of the former.

\lem 4
  Assume that $X$ is equidimensional and that $\dim P''(x_0)<r$.  Set
$T:=\Supp(\E/\N)$ and assume that $T\to Y$ is finite.  Then $\N$ is a
reduction of $\M$ if it is so over a dense Zariski open subset of $Y$.
 \pf
 Apply Corollary (10.7) of \cite{K-T, \p.225} as follows.  Let $A$ be
the local ring of $X$ at $x_0$.  Then $A$ is Noetherian, universally
catenary, and equidimensional.  Set $X_0:=\Specan A$.  Let $G$, $G'$,
and $G''$ be the quasi-coherent sheaves of graded algebras on $X_0$
associated to the stalks of $\S\E$, $\R\M$, and $\R\N$ at $x_0$.  Let
$P_0$, $P'_0$, and $P''_0$ be their ``Projan's.''  Set $M:=G$.  Then
$\Supp(\widetilde M)$ is equal to $P_0$, so it is equidimensional of
dimension $r_0$ with $r_0:=\dim X_0+p-1$.  Clearly, $\dim X_0=d+\dim
Y$.  Let $Y'$ be the closed subset of $X_0$ defined by the stalk of the
ideal of $T$.  Since $\N$ is equal to $\E$ off $T$, the stalks of
$\S\E$, $\R\M$, and $\R\N$ at $x_0$ become equal after localization
with respect to any analytic function on $X$ that vanishes along $T$;
hence, $G$, $G'$, and $G''$ are equal off $Y'$.  Moreover, since $T\to
Y$ is finite, the dimension of $Y'$ is bounded by that of $Y$.

Viewed as a $G''$-module, $M$ gives rise to a quasi-coherent sheaf on
$P''_0$.  The support $R$ of this sheaf is equal to $P''_0$; indeed, by
(6.4)(i) of \cite{K-T, \p.202}, $R$ is equal to the transform of $P_0$,
and by (2.6) of \cite{K-T, \p.183}, the latter is equal to $P''_0$.
Let $p''\:P''_0\to X_0$ be the structure map.  Then $p''^{-1}x_0$ is a
scheme, whose associated analytic space is $P''(x_0)$.  Hence
        $$\dim(p''^{-1}Y'\cap R)\le\dim p''^{-1}x_0 +\dim Y'
                \le r-1+\dim Y=r_0-1.$$
 (The corresponding bound in (10.7) of \cite{K-T, \p.225} is,
unfortunately, incorrectly stated because of a typographer's error;
however, the text suggests that the appropriate inequality is, in fact,
not strict.)  Moreover, if there is a component of dimension $r_0$ of
$p''^{-1}Y'\cap R$ (the $R$ is unnecessary in the present case), then
this component maps onto a component $Y'_1$ of $Y'$ such that $\dim
Y'_1=\dim Y$.

By hypothesis, $\R\M$ is a finitely generated module over $\R\N$
locally over a dense Zariski open subset $U$ of $Y$.  Let $Z$ be the
preimage in $T$ of $Y-U$, and $Z_0$ the closed subset of $X_0$ defined
by the stalk of the ideal of $Z$.  Then $\dim Z_0<\dim Y$, and so
$Y'_1$, if it exists, contains a point $\eta$ outside $Z_0$.  The
stalks $G'_{\eta}$ and $G''_{\eta}$ are localizations of the stalks of
$\R\M$ and $\R\N$ at $x_0$ with respect to a certain set of analytic
functions on $X$, including some that vanish on $Z$.  Hence $G'_{\eta}$
is a finitely generated module over $G''_{\eta}$; in the language of
\cite{K-T}, $G''_{\eta}$ is a reduction of $G'_{\eta}$ for $M_{\eta}$.
Finally, the preimage of $Y'$ in $P_0$ has no component of dimension
$r_0$, because $\dim Y'<\dim X_0$ and $\E$ is free of rank $p$.
Therefore, by (10.7) of \cite{K-T, \p.225}, $G'$ is a finitely
generated module over $G''$; in other words, at $x_0$ the stalk of
$\R\M$ is a finitely generated module over that of $\R\N$.  Hence
(after $X$ and $Y$ are replaced by neighborhoods of $x_0$ and $y_0$ if
necessary) $\N$ is a reduction of $\M$, and the proof is complete.

\prp 5
 Assume that $X$ is equidimensional.
 \part 1 If $y\mapsto e(y)$ vanishes, then  $\M=\E$ and $S=\emptyset$.
 \part 2 If $y\mapsto e(y)$ is constant on $Y$ and nonvanishing, then
$S\to Y$ is surjective, and $\M$ has a reduction generated by $r$
elements.
 \part 3 Assume that $Y$ is smooth, that $X/Y$ is flat with
Cohen--Macaulay fibers, and that $S\to Y$ is surjective.  If there
exists a reduction of $\M$ generated by $r$ elements, then, conversely,
$y\mapsto e(y)$ is constant on $Y$.
 \pf
 Consider (1).  Fix $y\in Y$.  Since $e(y)$ vanishes, a theorem of
Buchsbaum and Rim implies that the image of $\M$ in $\E(y)$ is all of
$\E(y)$.  Since $y$ is arbitrary, $\M=\E$ by Nakayama's lemma.  So
$S=\emptyset$ as $S:=\Supp(\E/\M)$.

Consider (2).  Fix $y\in Y$.  Since $e(y)$ doesn't vanish, the image of
$\M$ is, obviously, not all of $\E(y)$.  Hence $S\to Y$ is surjective.

After $X$ is replaced by a neighborhood of $x_0$ if necessary, there
exist $r$ elements of $\M$ whose images in $\E(y_0)$ generate a
reduction of the image of $\M$; see \cite{K-T, (6.6), \p.203} for
example.  Let $\N$ be the submodule of $\M$ generated by the elements.
Set $T:=\Supp(\E/\N)$.  Since $T(y_0)$ is finite, and since $(Y,y_0)$
is the germ of an analytic space, $T$ is finite over $Y$ after $X$ and
$Y$ are replaced by neighborhoods of $x_0$ and $y_0$ if necessary.
Hence, for every $y\in Y$, the Buchsbaum--Rim multiplicity $f(y)$ is
defined for the pair that $(\E,\M)$ induces on $X(y)$.  By
Proposition\Cs1), there is a dense Zariski open subset $U$ of $Y$ on
which $y\mapsto f(y)$ is constant and $f(y)\le f(y_0)$.  Then, for all
$y$ in $U$,
        $$e(y)\le f(y)\le f(y_0)=e(y_0)=e(y);$$
 the first relation holds because $\N\subseteq\M$, the second because
$y\in U$, the third by construction of $\N$, and the last by the
hypothesis.  Thus $e(y)=f(y)$.

Fix $y\in Y$.  For each $x\in X(y)$, let $e(x)$ denote the
Buchsbaum--Rim multiplicity at $x$ of the pair that $(\E,\M)$ induces
on $X(y)$, and let $f(x)$ denote that by $(\E,\N)$.  Of course, $e(x)$
vanishes if $x\notin S$, and $f(x)$ vanishes if $x\notin T$.  Clearly,
        $$e(y)=\sum_x e(x)\le \sum_x f(x)=f(y).$$
 Suppose $y\in U$.  Then $e(y)=f(y)$.  Hence $e(x)=f(x)$ for all $x\in
X(y)$.  Therefore, for all $x\in X(y)$, the stalk at $x$ of the image
of $\N$ in $\E(y)$ is a reduction of that of $\M$ by the generalized
theorem of Rees, Corollary (6.8)(a) of \cite{K-T, \p.207--8} with $M:=\S\E$.
Hence, for each $y$ in $U$, the image in $\E(y)$ of $\N$ is a reduction
of that of $\M$.  In particular,  $T=S$.

By Lemma\Cs2), there is a dense Zariski open subset of $Y$ over which
$\N$ is a reduction of $\M$.  On the other hand, $\dim P''(x_0)<r$
because $\N$ is generated by $r$ elements.  Finally, since $S\to Y$ is
surjective, so is $T\to Y$.  Therefore, Lemma\Cs4) implies that $\N$ is
a reduction of $\M$.

Consider (3).  Replacing $\M$ by its reduction, we may assume that $\M$
itself is generated by $r$ elements.  Let $\J$ denote the zeroth
Fitting ideal of $\E/\M$.  Then $\O_X/\J$ is supported by $S$.  Since
the codimension of $S$ is right, $\O_X/\J$ is Cohen--Macaulay.  Since
$Y$ is smooth and $S\to Y$ is finite and surjective, $\O_X/\J$ is
therefore flat over $Y$.  Hence the function
        $$y\mapsto \dim (F_*(\O_X/\J)(y))$$
 is constant on $Y$.  However, since the fibers $X(y)$ are
Cohen--Macaulay and since the formation of a Fitting ideal commutes
with base change,
        $$\dim (F_*(\O_X/\J)(y))=e(y)$$
 by some theorems of Buchsbaum and Rim \cite{B-R, 2.4 \p.207, 4.3 and
4.5 \p.223}.

\rmk 6
 Considerations involving the expression of  $e(y)$ as a sum of
intersection numbers suggest that, in Part (3) of Proposition\Cs5), the
Cohen--Macaulay hypothesis is unnecessary.  For example, it is
unnecessary when $p=1$ (so $\E=\O_X$  and $\M$ is an ideal); see
Theorem (2.2).

\lem 7
  The following conditions are equivalent:\smallskip
 \item i There exists a reduction of $\M$ generated by $r$ elements.
 \item ii The bound $\dim P'(x_0)<r$ obtains.
 \pf
 Indeed, assume (i) (replacing $X$ and $Y$ if necessary).  Say $\N$ is
the reduction of $\M$.  Then the inclusion $\R\N\into\R\M$ induces a
finite surjective map $P'\onto P''$.  Hence $\dim P'(x_0)=\dim
P''(x_0)$.  However, $\dim P''(x_0)<r$ because $\N$ is generated by $r$
elements.  Hence (ii) holds.

Conversely, assume (ii).  Then (i) follows, for example, from
\cite{K-T, (6.2)(iv), \p.201} applied with $\R(\M)$ for $G$ and for $M$
and with $\R(\N)$ for $G'$.  (The idea is simple.  Condition~(ii)
implies that there are $r$ hyperplanes in $P'$ whose intersection
misses the fiber $P'(x_0)$.  Let $\N$ be the submodule of $\M$
generated by the $r$ elements corresponding to these hyperplanes, and
$Z$ the subspace of $P'$ defined by the vanishing of these $r$
elements.  Then the central projection from $P'-Z$ to $P''$ restricts
to a finite map over $X-W$ where $W$ is the image of $Z$.  Hence, after
$X$ and $Y$ are replaced by neighborhoods of $x_0$ and $y_0$ if
necessary, $\N$ is a reduction of $\M$, and the proof is complete.)

\thm 8
 (Specialization of integral dependence)\enspace
 Assume that $X$ is equi\-di\-mensional, and that $y\mapsto e(y)$ is
constant on $Y$.  Let $h$ be a section of $\E$ whose image in $\E(y)$
is integrally dependent on the image of $\M$ for all $y$ in a dense
Zariski open subset of $Y$.  Then $h$ is integrally dependent on $\M$.
 \pf
 If $y\mapsto e(y)$ vanishes, then $\M=\E$ by Part (1) of
Proposition\Cs5), and so the assertion is trivial.  Assume $y\mapsto
e(y)$ is nonvanishing.  Then, by Part (2) of Proposition\Cs5), the map
$S\to Y$ is surjective, and (after $X$ is replaced by a neighborhood of
$x_0$ if necessary) there exists a reduction of $\M$ generated by $r$
elements.  So Lemma\Cs7) implies $\dim P'(x_0)<r$.

Let $\H$ be the submodule of $\E$ generated by $h$ and $\M$.  By
hypothesis, for all $y$ in a dense Zariski open subset of $Y$, the
image in $\E(y)$ of $\M$ is a reduction of the image of $\H$.  So, by
Lemma\Cs2), there is a smaller dense Zariski open subset of $Y$ over
which $\M$ is a reduction of $\H$.  Therefore, Lemma\Cs4) implies that
$\M$ is a reduction of $\H$, and the proof is complete.

\sctn The special fiber of the exceptional divisor

 Preserve the setup of Section (1).  Let $Z$ denote the analytic
subspace of $P$ defined by the sheaf of ideals in $\S\E$ generated by
$\M$.  Form the blowup $B$ of $P$ with respect to $Z$, and the
exceptional divisor $D$.  The main result of this section,
Theorem\Cs2), relates the condition $\dim D(y_0)<r$ to the constancy on
$Y$ of all the associated multiplicities $e^j(y)$ of the pair that $(\E,\M)$
induces on the fiber $X(y)$.  The definition of the $e^j(y)$ is
recalled below.  In particular, $e^0(y)$ is equal to $e(y)$, whose
constancy was studied in the last section, and part of that study will
be needed to prove Theorem\Cs2).  Not surprisingly, the constancy of
$e(y)$ alone does not imply the constancy of all the $e^j(y)$; one
instance where it doesn't is considered in Example\Cs3).  Half the
content of Theorem\Cs2) is provided by Lemma\Cs1), which gives a
geometric description of a dense (Zariski) open subset $U$ of $Y$ on
which all the $e^j(y)$ are locally constant.  In particular, Lemma\Cs1)
provides, in a second way, the open set $U$ needed in the proof of
Proposition~(1.1).

In the case $p=1$, there is only one possible nonzero associated
multiplicity, namely, $e^0(y)$.  In this case, $\M$ is the ideal on $X$
of $S$, and $B$ is the blowup of $X$ along $S$.  Theorem\Cs2) says
that, if $X$ is equimultiple along $S$, then the exceptional divisor is
equidimensional over $S$ (that is, if every fiber is empty or has the
minimal possible dimension $d-1$), and the converse holds if $Y$ is
smooth.  A version of the latter was proved in 1969 by Hironaka
(according to Remark~(2.6) in \cite{Lip, p.~121}); a few years later,
versions of the direct assertion were proved by Teissier (in \cite{T-1,
3.1, p.~327} and \cite{T1.5, I.1, p.~131, I.3, p.~133}) and by
Schikhoff (again according to \cite{Lip, p.~121}).

The main new technical ingredient in this section is intersection
theory.  Denote the first Chern classes of the tautological sheaves
$\O_{P'}(1)$ and $\O_{P}(1)$ by $\ell'$ and $\ell$.  Denote the blowup
of $P(y)$ with respect to $Z(y)$ by $B_y$, and denote the exceptional
divisor by $D_y$.  Finally, form the Segre numbers $s^i(y)$ of $Z(y)$
in $X(y)$:
  $$s^i(y):=\int\ell'^{i-1}\ell^{r-i}[D_y]\hbox{ for }i=1,\dots,r.$$
 Then $e^j(y)$ is defined as the sum of the first $r-j$ of the $s^i(y)$
in \cite{K-T, (7.1), \p.207}:
	$$e^j(y)=\sum_{i=1}^{r-j}s^i(y)\hbox{ for }j=0,\dots,r-1.$$
 (In fact, in \cite{K-T, (7.1)}, the sum starts with an additional term
$s^0(y)$, but that term clearly vanishes here.)  With this definition,
	$$e(y)=e^0(y)$$
 because of \cite{K-T, (5.1), \p.191} and \cite{K-T, (5.7), \p.207}.
The projection formula with respect to the map $D_y\to Z(y)$ yields
	$$s^i(y)=0 \hbox{ for } i<d;$$
 in particular, $e^{r-d}(y)=s^d(y)$.
 All the $e^j(y)$ are constant if and only if all the $s^i(y)$ are.
However, the $e^j(y)$ are upper semicontinuous, whereas the $s^i(y)$
needn't be; see Example\Cs3).

\lem 1
 Let $U$ be the open subset of $y$ in $Y$ such that $\dim D(y)<r$ and
$Y$ is smooth at $y$.  Then on $U$ all the functions $y\mapsto e^j(y)$
and $y\mapsto s^i(y)$ are locally constant.
 \pf
 For each $y\in Y$, the fiber $B(y)$ contains the blowup $B_y$ as a
closed subscheme, and the intersection $D\cap B_y$ is equal to the
exceptional divisor $D_y$.  Hence the intersection product $D\cdot
[B_y]$ is equal to the fundamental cycle $[D_y]$.  Now, $B(y)-D(y)$ is
equal to $B_y-D_y$.  So, if $\dim D(y)<r$, then $[B(y)]$ is equal to
$[B_y]$, and $D\cdot [B_y]$ is equal to $[D(y)]$; hence
  $$s^i(y)=\int\ell'^{i-1}\ell^{r-i}[D(y)]\hbox{ for }i\ge1.$$
 If $y\in U$, then the embedding $\iota_y\:y\into Y$ is regular; so the
operation of pullback along $\iota_y$ commutes with that of pushforth
along the proper map $D\to Y$.  Hence, $y\mapsto s^i(y)$ is constant on
each connected component of $U$ for $i\ge1$; compare with \cite{Fulton,
10.2, \p.180}.  The proof is now complete.

\thm 2
  Assume that $X$ is equidimensional.  If the function $y\mapsto
e^j(y)$ is constant on $Y$ for $0\le j<p$, then the central fiber
$D(y_0)$ of the exceptional divisor $D$ of the blowup $B$ of $P$ is
empty or has the minimal possible dimension, $r-1$.  Furthermore, the
converse holds if, in addition, $Y$ is smooth at $y_0$.
 \pf
 The converse follows immediately from Lemma\Cs1), applied after $Y$ is
replaced by $U$.  So assume that $e^j(y)$ is constant on $Y$ for all
$j$.  By way of contradiction, suppose that $D(y_0)$ has a component
$D'(y_0)$ of dimension $r$ or more.  Replacing $X$ by a neighborhood of
$x_0$ if necessary, we may assume that $x_0$ is the unique point of
$S(y_0)$.  Then $D(y_0)=D(x_0)$.  By Proposition~(1.5)(2), after $X$
and $Y$ are replaced by neighborhoods of $x_0$ and $y_0$ if necessary,
there exists a reduction $\N$ of $\M$ generated by $r$ elements.  So
Lemma~(1.7) yields $\dim P'(x_0)<r$.  However, $P'(x_0)$ contains the
image of $D'(x_0)$ under the projection of $B$ onto $P'$.  Hence the
fibers of the map $D'(x_0)\to P'(x_0)$ all have dimension at least 1.
However, these fibers are embedded in $P(x_0)$ by the blowup map $B\to
P$ for the following reason: by definition of $Z$, its ideal sheaf is a
quotient of the pullback of $\M$ to $P$, and so the Rees algebra of the
ideal sheaf is a quotient of $\R\M$; correspondingly, $B$ is embedded
in $P'\x P$, and the second projection restricts to the blowup map
$B\to P$.  Now, if $p=1$, then $P=X$, and so $P(x_0)$ can contain no
subspace of dimension at least 1.  Thus, if $p=1$, then the assertion
holds.

The proof proceeds by induction on $p$.  Suppose $p>1$.  Let $g$ be a
general section of $\E:=\O_X^p$.  Set $\E':=\E/g$ and let $\M'$ denote
the image of $\M$.  Then $\E'$ is free of rank $p-1$.  Set
$Q:=\Projan(\S(\E'))$ and $Z':=Z\cap Q$.  Obviously $Z'$ is defined by
the sheaf of ideals generated by $\M'$.  Since $g$ is general, the
preimage of $Q$ in $B$ is equal to the blowup $C$ of $Q$ along $Z\cap
Q$, and $D\cap C$ is the exceptional divisor $E$, at least after $Y$ is
replaced by a neighborhood of $y_0$.
 Since the fibers of the map $D'(x_0)\to P'(x_0) $ have dimension at
least 1 and are embedded in $P$ by the blowup map, it follows that
$Q$ must intersect these fibers.
 Hence the fiber $E(y_0)$ has dimension $r-1$ or more.  For
convenience, denote the $j$th associated multiplicity of the pair that
$(\E',\M')$ induces on $X(y)$ by $e'^j(y)$.  Since $g$ is general, it
follows from \cite{K-T, (7.1) and (7.2)(iv), \pp.207--8} that $e'^j(y)$
is equal to $e^{j+1}(y)$ for $0\le j\le r-1$ and for all $y$ in a
(Zariski) open neighborhood of $y_0$; replace $Y$ by this neighborhood.
By hypothesis, $e^{j+1}(y)$ is constant on $Y$.  Hence $y\mapsto
e'^j(y)$ is constant.  Thus the induction hypothesis is contradicted,
and so $D'(y_0)$ does not exist.  The proof is now complete.

\eg 3
 Consider the example of Henry and Merle \cite{H&M83, \p.578--9}.  In it,
the parameter space $Y$ is the affine line ${\bf C}$, and the total
space $X$ is cut out of ${\bf C}^4\times Y$ by two equations,
	$$X:X_1^2+X_2^2+X_3^2+yX_4=0,X_1^4+X_2^4+X_3^4+X_4^2=0.$$
 Set $p:=2$, and let $\M$ be the Jacobian module, the column space of
the Jacobian matrix with respect to the $X_i$.  Henry and Merle proved
that $X$ is Whitney equisingular along the $Y$-axis at the origin.
Hence, by \cite{G-4, 1.3, p.~211, 2.6, p.~215}, the function $y\mapsto
e(y)$ is constant on $Y$.  In fact, in the case at hand, it is not hard
to see via a direct computation that $e(y)=36$ for all $y$.

On the other hand, $e^1(y)=0$ if $y\ne0$ because, obviously, $X(y)$ has
embedding dimension 3 at the origin.  However, $e^1(0)\ne0$ because
$e^1(0)$ is the multiplicity of the ideal obtained by taking a generic
linear combination of the rows of the Jacobian matrix of $X(0)$; in
fact, it is easy to see that $e^1(0)=4$.  Hence, by Theorem (2.2),
$D(0)$ must have a ``vertical'' component.  In fact, it is also not
hard to see that $D(0)$ consists of two components, one of which maps
onto the fiber of $P$ over the origin in $X(0)$.

Finally, $e(y)$ is equal to $e^1(y)+s^3(y)$; so this sum is constant.
On the other hand, $e^1(y)$ is upper semicontinuous.  Therefore,
$s^3(y)$ is lower semicontinuous.

\sctn Integral dependence and strict dependence

In the next sections, we'll study Whitney's Condition A, Thom's
Condition A$_f$, and Henry, Merle and Sabbah's Condition \Wf, which
concern limiting tangent hyperplanes at a singular point of a complex
analytic space.  To prepare further for this study, in this section and
in part of the next one, we'll recall and develop some material from
\cite{G-2} and \cite{G-1}.  In \cite{T-1}, Teissier made a similar
study in the case of families of hypersurfaces with isolated
singularities, and his work has been a model for ours.

Let $(X,0)$ be the germ of a complex analytic space, and $\E:=\O_X^p$
a free module of rank $p$ at least 1. Let $\M$ be a coherent submodule
of $\E$, and $h$ a section of $\E$.  Given a map of germs
$\vf\:(\IC,0)\to(X,0)$, denote by $h\o\vf$ the induced section of the
pullback $\vf^*\E$, or $\O_{\IC}^p$, and by $\M\o\vf$ the induced
submodule.  Call $h$ {\it integrally dependent\/} (resp., {\it strictly
dependent\/}) on $\M$ at $0$ if, for every $\vf$, the section $h\o\vf$
of $\vf^*\E$ is a section of $\M\o\vf$ (resp., of $m_1(\M\o\vf)$, where
$m_1$ is the maximal ideal of $0$ in $\IC$).  The submodule of $\E$
generated by all such $h$ will be denoted by $\?{\M}$, resp., by $\sdep
\M$ (the notation `$\sdep \M$' is a change from \cite{G-1}).

To check for integral (resp., strict) dependence, it suffices to use
only those $\vf$ whose image meets any given dense Zariski open subset
of $X$.  Indeed, if $h\o\vf$ is not a section of $\M\o\vf$ (resp., of
$m_1(\M\o\vf)$, then $\vf$ can be tweaked, preserving this condition, so
that the image of $\vf$ does meet the given open set (see the proof of
Prop.~1.7 on \p.304 in \cite{G-2}).

Let $\N$ be a coherent submodule of $\M$.  Then $\M\subset\?{\N}$ if and
only if $\N$ is a reduction of $\M$ in the sense of Section~1 (after $X$
is replaced by a neighborhood of 0 if necessary).  Indeed, the present
definition of integral closure is taken from \cite{G-2, 1.3, \p.303}.
This definition is shown, on the middle of \p.305 in \cite{G-2}, to be
equivalent to Rees's definition \cite{R87, \p.435}.  Hence,
Theorem~1$\cdot$5 in \cite{R87, \p.437} yields the assertion.

The following result is a simple, but useful, observation.

\prp1
 If $\N\subset\M\subset\?{\N}$, then $\?\M=\?{\N}$ and
$\sdep{\M}=\sdep{\N}$.
 \pf
 For any map  $\vf\:(\IC,0)\to(X,0)$, the hypothesis yields
	$$\N\o\vf\subset\M\o\vf\subset\?\N\o\vf.$$
 By definition, the third term is equal to the first.  Hence the first
term is equal to the second.  Therefore, the definitions yield the
assertions.
	\medbreak
The next result gives one useful connection between the notions of
integral dependence and strict dependence.

\prp2
 Fix a set $C$ of generators of $\M$.  Then the following conditions on
$\M$ and $\N$ are equivalent:
 \smallskip
 \item i  $\N\subset\sdep{\M}$;
 \item ii  $\M\subset\?{\N'}$ for every coherent submodule $\N'$ of
$\M$ such that $\N+\N'=\M$;
 \item iii $\M\subset\?{\N'}$ for every submodule $\N'$ of $\M$ such
that $\N'$ is generated by a subset of $C$ and $\N+\N'=\M$.
 \pf
 Assume (i).  To prove (ii), take any map $\vf\:(\IC,0)\to(X,0)$.  Then
	$$\M\o\vf=\N\o\vf+{\N'}\o\vf\subset m_1(\M\o\vf)+\N'\o\vf.$$
 By Nakayama's lemma, $\M\o\vf=\N'\o\vf$.  So (ii) holds.  Trivially,
(ii) implies (iii).

Finally, assume (i) fails.  Then there exists a $\vf$ such that
$\N\o\vf$ is not contained in $m_1(\M\o\vf)$.  Let $h$ be a section of
$\N$ such that $h\o\vf$ is not contained in $m_1(\M\o\vf)$.  Supplement
$h$ by elements of $C$ to obtain a basis of the vector space
$\M\o\vf\big/m_1(\M\o\vf)$, or what is the same, a basis of $\M/m\M$
where $m$ is the maximal ideal of $0$ in $\O_X$.  Let $\N'$ be the
submodule of $\M$ generated by these elements of $C$.  Then, by
Nakayama's lemma, $\N+\N'=\M$.  Moreover, by construction, $h\o\vf$ is
not contained in $\N'\o\vf$.  Thus (iii) fails, and the proof is
complete. \medbreak

The following lemma is a useful generalization of Proposition~1.6 in
\cite{G-1}, and the following proof is a little different.

\lem3 For a section $h$ of $\E:=\O_X^p$ to be integrally dependent (resp.,
strictly dependent) on $\M$ at $0$, it is necessary that, for all maps
$\vf\:(\IC,0)\to(X,0)$ and $\psi\:(\IC,0)\to (\Hom(\IC^p,\IC),\lambda)$
with $\lambda\ne0$, the function $\psi(h\o\vf)$ on $\IC$ belong to the
ideal $\psi({\M}\o\vf)$ (resp., to {$m_1\psi(\M\o\vf)$}).

  Conversely, it is sufficient that this condition obtain for every
$\vf$ whose image meets any given dense Zariski open subset of $X$.
Furthermore, if $0$ lies in the cosupport $\Supp(\E/\M)$ of $\M$, then
it is sufficient that the condition obtain for every such $\vf$ and for
every $\psi$ that carries $\M\o\vf$ into $m_1$.
 \pf
 The first assertion follows directly from the definitions.

Conversely, given any $\vf$, (after $X$ is replaced by a neighborhood of
0 if necessary) there exists a basis $e_1,\dots e_p$ for $\vf^*\E$ such
that $\M\o\vf$ is equal to the submodule generated by
$t^{n_1}e_1,\dots,t^{n_r}e_r$ for suitable integers $n_i$ and $r$, where
$t$ is the coordinate function on $\IC$.  Say $h\o\vf$ expands as
$a_1e_1+\cdots+a_pe_p$.  Then clearly $h\o\vf$ is a section of $\M\o\vf$
(resp., of $m_1(\M\o\vf)$) if and only if $a_i=b_it^{n_i}$ (resp.,
$a_i=b_it^{n_i+1}$) for a suitable $b_i$ for $1\le i\le r$ and $a_i=0$
for $r<i\le p$.  Form the dual basis $e_1',\dots,e_p'$.  Then
$a_i=e_i'(h\o\vf)$ for all $i$.  Moreover, $e_i'(\M\o\vf)$ is equal to
the ideal generated by $t^{n_i}$ for $1\le i\le r$, and to 0 for $r<i\le
p$.

Hence $h\o\vf$ lies in $\M\o\vf$ (resp., in $m_1(\M\o\vf)$) if (and only
if) the condition obtains for the $p$ maps $\psi$ corresponding to
$e_1',\dots,e_p'$, and for each of these $\psi$, obviously
$\lambda\ne0$.  Thus the second assertion holds.

Suppose $0\in\Supp(\E/\M)$.  Then either $r=p$ and $n_j>0$ for some $j$,
or $r<p$.  Suppose first $r<p$.  Fix $i$ with $r<i\le p$, and let
$\psi$ correspond to $e_i'$.  Then the condition implies that $a_i=0$.
Now, fix $i$ with $1\le i\le r$, and let $\psi$ correspond to
$te_i'+e_{r+1}'$.  Then $\lambda\ne0$, and $\psi$  carries $\M\o\vf$
into $m_1$.  The condition implies that $ta_i=b_it^{{n_i}+1}$ (resp.,
$ta_i=b_it^{n_i+2}$) for some $b_i$.  Thus the third assertion holds
when $r<p$.

Suppose $r=p$.  Reorder the $e_i$ so that $n_i\le n_{i+1}$ for each $i$.
Say $n_j=0$, but $n_{j+1}>0$.  Fix $i$ with $j<i\le p$, and let $\psi$
correspond to $e_i'$.  Then $\lambda\ne0$, and $\psi$ carries
$\M\o\vf$ into $m_1$.  The condition implies that $a_i=b_it^{n_i}$
(resp., $a_i=b_it^{n_i+1}$) for some $b_i$.  Now, fix $i$ with $1\le
i\le j$.  Then $n_i=0$, so $a_i=a_it^{n_i}$.  Thus $h$ is integrally
dependent on $\M$ at $0$.

To handle strict dependence, let $\psi$ correspond to $te_i'+e_{j+1}'$.
Then $\lambda\ne0$, and $\psi$  carries $\M\o\vf$ into $m_1$.  The
condition implies that $ta_i+a_{j+1}=b_i(t+t^{n_{j+1}})$ for some
$b_i$.  Now, $n_{j+1}>0$ and $a_{j+1}=b_{j+1}t^{n_{j+1}+1}$.  Hence,
$a_i=b_i't^{n_i+1}$ for a suitable $b_i'$.  Thus the third assertion
holds, and the proof is complete.
 \medbreak

It is often convenient to work on the space $P$ of Section~1;
obviously, $P=X\x\IP^{p-1}$ since $\E:=\O_X^p$.  The section $h$ of
$\E$ and the submodule $\M$ of $\E$ generate ideals on $P$; denote them
by $\rho(h)$ and $\rho(\M)$.  Note that $\rho(h)$ is locally principal.
The next result gives a translation of the two notions of dependence
into this context, thereby reducing the study of dependence on the
module $\M$ on the germ $(X,0)$ to that of the ideal $\rho(\M)$ on the
more global space $P$.

\prp4
 A necessary and sufficient condition for a section $h$ of $\E$ to be
integrally dependent (resp., strictly dependent) on $\M$ at $0$ is that,
at each point of $\IV(\rho(\M))$ lying over $0\in X$, a generator of
$\rho(h)$ be integrally dependent (resp., strictly dependent) on
$\rho(\M)$.
 \pf
 To give a map $\phi\:(\IC,0)\to(P,(0,l))$ is the same as to give a
pair of maps $\vf\:(\IC,0)\to(X,0)$ and $\psi\:(\IC,0)\to
(\Hom(\IC^p,\IC),\lambda)$ where $\lambda$ corresponds to $l$ (although
$\psi$ is determined only up to multiplication by a function that
doesn't vanish at $0\in\IC$).  Hence, the assertion follows from
Lemma~\Cs3).
 \medbreak

It is also convenient to work with the normalized blowup, with its
structure map,
	$$\pi\:NB_{\rho(\M)}(P)\to P,$$
 and with its exceptional divisor $E$.  (After replacing $X$ by a
neighborhood of 0 if necessary, we may assume that each component of $E$
meets the fiber over $0$.)  The next result relates the two notions of
dependence to vanishing of the ideal $\rho(h)\o\pi$ on the components of
$E$.

\prp5 Let $h$ be a section of $\E$, and $Y$ a closed analytic subset of
the image of $E$ in $X$.
 \part1 A necessary and sufficient condition for $h$ to be integrally
dependent on $\M$ at $0$ is that, along each component of $E$, the ideal
$\rho(h)\o\pi$ vanish to order at least the order of vanishing of
$\rho(\M)\o\pi$.
 \part2 A necessary and sufficient condition for $h$ to be strictly
dependent on $\M$ at every $y\in Y$ is that, along each component $V$ of
$E$, the ideal $\rho(h)\o\pi$ lie in the product $\I(Y,V)\rho(\M)\o\pi$,
where $\I(Y,V)$ denotes the ideal of the reduced preimage of $Y$ in $V$;
in particular, if $V$ projects into $Y$, then this condition simply
requires the ideal $\rho(h)\o\pi$ to vanish to order strictly greater
than the order of vanishing of $\rho(\M)\o\pi$.
 \pf
 Consider (1).  Proposition~\Cs4) reduces the assertion to the case of
an ideal, and this case is treated in \cite{T-2,\p.330, 1.4 Prop.~ 2}.

 Consider (2).  At each  $b\in V$, the ideal $\rho(\M)\o\pi$ is
generated by a single section $g\o\pi$ where $g$ is a suitable section
of $\M$, and the ideal $\rho(h)\o\pi$ is generated by a multiple
$k(\rho(g)\o\pi)$ where $k$ is a meromorphic function.  In these terms,
the condition in (2) says that $k$ is holomorphic and vanishes at $b$ if
$b$ projects into $Y$.

Hence the condition in (2) holds if and only if, for every map
 $$\beta\:(\IC,0)\to (NB_{\rho(\M)}(P),b)$$
 such that $\phi:=\pi\o\beta$ is not constant and such that the image of
$\beta$ meets the complement of $\IV(\rho(\M))$, the function
$k\o\beta$ vanishes at $0\in\IC$.  Now, $\rho(h)\o\phi$ is generated at
$0$ by $(k\o\beta)(g\o\phi)$ if $k$ is holomorphic; moreover, $k$ is
holomorphic, if $h$ is integrally dependent on $\M$ at image of $b$ by
(1).  Furthermore, we can factor any map $\phi\:(\IC,0)\to (P,(y,l))$
through $NB_{\rho(\M)}(P)$.  Therefore, the assertion follows from
Proposition~\Cs4).

\sctn Whitney's Condition \A

In this section, we use the theory developed in the preceding
sections to study Whitney's Condition A.  After introducing the setup,
we prove a lemma, which relates limit tangent hyperplanes with the
notions of strict dependence and integral dependence; the statement and
proof are, more or less, found in Section 2 of \cite{G-1}.  Then we
prove the main result of the section, Theorem~\Cs2), which asserts that
Whitney's Condition A holds under the constancy of certain
Buchsbaum--Rim multiplicities on the fibers.  Finally, we illustrate
the theorem with two examples.

Let $(X,0)$ be a complex analytic subgerm of $(\IC^n,0)$ defined by the
vanishing of a map of germs $F\:(\IC^n,0)\to (\IC^p,0)$.  Call the
$\O_X$-submodule of the normal module to $X$ in $\IC^n$ generated by
all the partial derivatives of $F$ the ({\it absolute\/}) {\it Jacobian
module} of $F$, and denote it by $JM(F)$; more precisely, $JM(F)$ is
the image of the canonical map,
        $$\IHom_X(\Omega^1_{\IC^n}|X,\O_X)\to\IHom_X(\I/\I^2,\O_X),$$
where $\I$ is the ideal of $X$ in $\IC^n$.  Since $\I$ is generated by
the $p$ coordinate functions of $F$, the displayed map is given by the
Jacobian matrix $DF$, and $JM(F)$ is simply the submodule of the free
module $\O_X^p$ generated by the columns of $DF$.  Note in passing that this
module $\O_X^p$ contains the target
$\IHom_X(\I/\I^2,\O_X)$, which is an abstract $\O_X$-module and is known
as the {\it normal module\/} of $X$ in $\IC^n$, but the embedding
depends on the choice of the $p$ generators of $\I$; moreover, this
embedding is an isomorphism if $X$ is a complete intersection of
codimension $p$, but not in general.

Given an analytic map germ $g\:(\IC^n,0)\to (\IC^l,0)$, let
$JM(F)_g$ denote the submodule of $JM(F)$ generated by the ``partials''
$\pd Fv$ for all vector fields $v$ on $\IC^n$ tangent to the fibers of
$g$, that is, for all $v$ that map to the 0-field on $\IC^l$; call
$JM(F)_g$ the {\it relative Jacobian module\/} with respect to $g$.
For example, if $g$ is the projection onto the space of the last $l$
variables of $\IC^n$, then $JM(F)_g$ is simply the submodule generated
by all the partial derivatives of $F$ with respect to the first $n-l$
variables.

Call a hyperplane in $\IC^n$ through $0$ a {\it limit tangent
hyperplane\/} of $(X,0)$ if it is the limit of hyperplanes tangent to
$X$ at nonsingular points approaching 0 along an analytic arc.  Now,
let $(S,0)$ be a smooth subgerm of $(\IC^n,0)$ defined by the vanishing
of a map of germs $g\:(\IC^n,0)\to (\IC^l,0)$ with $l=n-m$ where
$m:=\dim S$, and let $T_0S$ denote its tangent space at $0$.  Finally,
denote the singular locus of $X$ by $\Sigma$.

The following lemma describes the limit tangent hyperplanes in general
and those that contain $T_0S$ in particular in terms of Jacobian
modules.  The lemma and its proof are, more or less, the statement and
proof of Theorem~2.4 of \cite{G-1}.

\lem1 Preserve the conditions above.
 \part1  A hyperplane $H$, defined by the vanishing of a linear
function $h\:\IC^n\to \IC$, is a limit tangent hyperplane of $(X,0)$ if
and only if $JM(F)_h$ is not a reduction of $JM(F)$.
 \part2 No hyperplane containing $T_0S$ is a limit tangent hyperplane
of $(X,0)$ if $JM(F)_g$ is a reduction of $JM(F)$.
 \part3 Every limit tangent hyperplane of $(X,0)$ contains $T_0S$~--- in
other words, the pair $(X-\Sigma,S)$ satisfies {\rm Whitney's Condition
A} at the origin~--- if and only if $JM(F)_g$ is contained in
$\sdep{JM(F)}$.
 \pf
 A hyperplane $M$ is a limit tangent hyperplane of $(X,0)$ if and
only if there exist maps $\vf(t)\:(\IC,0)\to(X,0)$ and
$\psi(t)\:(\IC,0)\to (\Hom(\IC^p,\IC),\lambda)$ with $\lambda\ne0$ such
that $\vf(t)$ is a nonsingular point of $X$ for $t\neq0$ and such that,
for a suitable $k$, the limit,
        $$\lim_{t\to0}(1/t^k)\bigl(\psi(t)DF(\vf(t))\bigr),$$
 exists and is a conormal vector to $M$.  This condition means that,
given a vector field $v$ on $\IC^n$, the vector $v(0)$ lies in $M$ if
and only if
 $$\psi(t)DF(\vf(t))v(\vf(t))\in m_1\bigl(\psi(t)JM(F)\o\vf(t)\bigr),$$
 where $m_1$ is the maximal ideal of $0$ in $\IC$.  Obviously,
$DF(\vf(t))v(\vf(t))$ is equal to $(\pd Fv)\o\vf(t)$.  Hence, $v(0)\in
M$ if and only if
        $$\psi(t)(\pd Fv)\o\vf(t)\in m_1(\psi(t)JM(F)\o\vf(t)).\dno1.1$$

Consider (1).  If the vector field $v$ is tangent to the fibers of $h$,
then $v(0)\in H$, and every vector in $H$ is a $v(0)$
for some such $v$.  Hence, if $H=M$, then \Cs1.1) holds for every $v$
tangent to the fibers of $h$; whence, $JM(F)_h$ is not a reduction of
$JM(F)$, thanks to Lemma~(3.3) applied with $\pd Fv$ for $h$ and with
$JM(F)$ for $\M$.  Conversely, if $JM(F)_h$ is not a reduction of
$JM(F)$, then, by Lemma~(3.3), there exists a pair of maps $\vf(t)$ and
$\psi(t)$ such that \Cs1.1) holds for every $v$ tangent to the fibers
of $h$; whence, the corresponding  $M$ contains every vector in $H$,
and so $M=H$.  Thus (1) holds.

Consider (2) and (3).  Since $S$ is smooth and $l=n-m$, the germ $g$ is
a submersion.  Hence, if the vector field $v$ is tangent to the fibers
of $g$, then $v(0)\in T_0S$, and every vector in $T_0S$ is a $v(0)$ for
some such $v$.  If $JM(F)_g$ is a reduction of $JM(F)$, then \Cs 1.1)
fails for some such $v$ by Lemma~(3.3), and therefore $v(0)\notin M$;
thus (2) holds.  Finally, \Cs 1.1) implies that $v(0)$ lies in every
limit tangent $M$ if and only $\pd Fv$ is contained in $\sdep{JM(F)}$,
thanks to Lemma~(3.3); thus (3) holds.  The proof of the lemma is now
complete.
 \medbreak

The following theorem gives a sufficient fiberwise numerical criterion
for the condition to hold.  The proof involves a delicate interplay
among the absolute and two relative Jacobian modules of $F$.
  \thm2
 Let $Y:=\IC^m$ be the space of the first $m$ coordinates in $\IC^n$
where $1\le m<n$, and set $l:=n-m$.  Assume that, under the projection
$r\:\IC^n\to Y$, the subspace $X$ of $\IC^n$ becomes the total space of
a family of complete intersections $X(y)$ of codimension $p$ defined by
the maps $F(y)\:\IC^l\to\IC^p$ given by $F(y)(z):=F(y,z)$.  Assume that
the $X(y)$ have isolated singularities, which trace out the smooth
subgerm $(S,0)$ of $(X,0)$.  Finally, let $e^j(y)$ be the $j$th
associated Buchsbaum--Rim multiplicity of the Jacobian module $JM(F(y))$
in $\O_{X(y)}^p$, and assume that the function $y\mapsto e^j(y)$ is
constant on $(Y,0)$ for $0\le j< p$.  Then $(X-S,S)$ satisfies Whitney's
Condition \A\ along $S$.
 \pf
 Form the ideal $\rho(JM(F)_r)$ on $X\x\IP^{p-1}$, form the
corresponding normalized blowup, and form its exceptional divisor $E$.
Then every component of $E$ projects onto $Y$; indeed, this conclusion
follows from Theorem~(2.2) because, by hypothesis, the functions
$y\mapsto e^j(y)$ are constant.

Since $S$ is smooth of dimension $m$, there is a map germ
$g\:(\IC^n,0)\to (\IC^l,0)$ such that $(S,0)=(g^{-1}0,0)$.  Moreover,
by the generic Whitney lemma, Whitney's Condition A holds on a dense Zariski
open subset  $U$ of $S$.  Hence, for  $s\in U$,
        $$JM(F)_g\subset\sdep{JM(F)}\hbox{ at }s\dno2.1$$
 by Proposition~\Cs1).  Replacing $U$ by a smaller subset, we may assume
that the map $S\to Y$ is unramified at $s$.  Then, at $s$, the sum
$JM(F)_g+JM(F)_r$ is all of $JM(F)$.  Hence, at $s$,
        $$JM(F)\subset\?{JM(F)_r}$$
 by (i)$\Rightarrow$(ii) of Proposition~(3.2).

Hence, this inclusion holds everywhere on $X$ by Proposition~(3.5)(1)
since every component of $E$ projects onto $Y$; apply the proposition
twice, first the necessity assertion with $U$ for $Y$, and then the
sufficiency assertion.  (In fact, here we could appeal to Theorem~(1.8)
instead, and thus use only the constancy of $e^0(y)$; however, the
constancy of all the $e^j(y)$ is used in an essential way in the next
paragraph.)  Therefore, by Proposition~(3.1),
        $$\sdep{JM(F)}=\sdep{JM(F)_r}\dno2.2$$
 everywhere on $X$.

Again, since every component of $E$ projects onto $Y$,
Proposition~(3.5)(2) implies that the inclusion,
        $$JM(F)_g\subset\sdep{JM(F)_r},$$
 holds everywhere on $S$ because it holds on $U$ by virtue of \Cs2.1)
and \Cs2.2).  Therefore, again by virtue of \Cs2.2), the inclusion
\Cs2.1) holds along $S$.  Consequently, Proposition~\Cs1) implies that
the pair $(X-S,S)$ satisfies Whitney's Condition A along $S$, and the
proof is complete.

\eg3 The fiberwise numerical criterion of Theorem~\Cs2), although
sufficient, is not necessary.  In fact, it is impossible to have a
necessary and sufficient numerical criterion for Whitney's Condition A
that depends only on the members of the family.  This observation was
made by Trotman \cite{Trot, Prop.~5.1, p.~147} on the basis of the
following example of his:
        $$X:w^a=y^bv^c+v^d$$
 and $S$ is the $y$-axis.  Here, different values of $b$ can give
essentially different parameterizations of the same collection of plane
curves.  Trotman determined when Whitney's Condition A is satisfied,
and when it isn't.

Let's look, from our point of view, at a special case of Trotman's
example, the ``Whitney umbrella of type $b$,''
        $$X:w^2-v^3+v^2y^b=0.$$
 Since $X$ is the total space of a one-parameter family of plane
curves, $m=1$ and $l=2$ and $p=1$; furthermore, $e^0(y)$ is simply the
(ordinary) multiplicity of the Jacobian ideal.  Here, $X(y)$ is a nodal
cubic for $y\neq0$, and $X(0)$ is a cuspidal cubic.  So $e^0(0)=3$ and
$e^0(y)=2$ for $y\neq0$, as is easy to check.  Finally, it is well
known that Whitney's Condition A is satisfied if $b\ge2$, but not if
$b=1$.  This fact will now be checked as an illustration of the use of
Lemma~\Cs1)(3).

Set $F:=w^2-v^3+v^2y^b$.  Let $g\:(\IC^3,0)\to (\IC^2,0)$ be the
projection onto the $(w,v)$-plane.  Then $JM(F)_g$ is generated by the
partial derivative $\pd Fy$, and so we have to show that $\pd Fy\in
\sdep{JM(F)}$ holds if $b\ge2$, but fails if $b=1$.  So consider a map
$\vf(t)\:(\IC,0)\to(X,0)$, say with coordinate functions,
        $$w(t)=\alpha t^i+\cdots,\ v(t)=\beta t^j+\cdots,\
                y(t)=\gamma t^k+\cdots,$$
 with $i,j,k>0$.  Since $DF=(2w,\ v(3v+2y^b),\ bv^2y^{b-1})$, we're
asking about the condition,
        $$2j+(b-1)k > \min(i,\ 2j+j'),\dno2.3$$
  where $j'\ge0$ and $j'=0$ unless $j=bk$ and $3\beta+2\gamma^b=0$.  We
have to establish this condition if $b\ge2$, and show that it fails if
$b=1$ for a suitable choice of $\vf(t)$. Now,
        $$\alpha^2 t^{2i}+\cdots=(\beta^2 t^2j+\cdots)
                (\beta t^j-\gamma^b t^{bk} +\cdots). \dno2.4$$
 So there are three cases to consider.  First, suppose that $j>bk$.
Then \Cs2.4) implies that $2i=2j+bk<3j$.  So $i<3j/2$.  So \Cs2.3)
holds for any $b\ge0$.  Second, suppose either that $j<bk$ or that
$j=bk$ and $\beta\neq\gamma^b$.  Then \Cs2.4) implies that $2i=3j$.  So
$i=3j/2$.  So again \Cs2.3) holds for any $b\ge0$.  Finally, suppose
that $j=bk$ and $\beta=\gamma^b$.  If $b\ge2$, then \Cs2.3) holds
because $k>0$ and $j'=0$.  However, if $b=1$, then \Cs2.3) need not
hold.  For instance, \Cs2.3) does not hold if
        $$w(t)=t^2,\ v(t)=t,\ y(t)=t-t^2,$$
 although \Cs2.4) does hold.  Thus Whitney's Condition A is satisfied
if $b\ge2$, but not if $b=1$.

\eg4 In Theorem~\Cs2), the smooth subgerm $(S,0)$ of $(X,0)$ is
customarily taken to be the plane of the first $m$ variables, but it
needn't be.  In fact, the projection of $(S,0)$ onto $(Y,0)$ may be
allowed to ramify at $0$.  For example, suppose that $X$ is defined as
follows:
        $$X:(w^2-y)^2-x^2=0.$$
 Then $X(y)$ is a binodal quartic for $y\neq0$ and $X(0)$ is a tacnodal
quartic.  Moreover, $S$ is the parabola,
        $$S:w^2=y,\ x=0.$$
 It is easy to check that $e(y)=4$ for all $y$.  Hence, by
Theorem~\Cs2), Whitney's Condition A is satisfied.  In fact, it is
obvious geometrically that the condition is satisfied.  Indeed,
$X=X_+\cup X_-$ where
        $$X_+:(w^2-y)+x=0\and X_-:(w^2-y)-x=0.$$
 The two components $X_+$ and $X_-$ are smooth, and they meet
transversally along their intersection, which is $S$.  Therefore, each
limit tangent hyperplane at 0 is either the tangent plane to $X_+$ or
that to $X_-$, so contains the tangent line to $S$.

\sctn Thom's Condition \Af

In this section, we use the theory developed in Sections 2 to 4 to
study Thom's Condition \Af.  After introducing the setup, we prove a
lemma, which is similar to Lemma~(4.1), and relates limit tangent
hyperplanes to level hypersurfaces with the notions of strict
dependence and integral dependence.  Then we prove a generalization of
the L\^e--Saito theorem.  Finally, we prove some variations of a
special case of a recent result of Brian\c con, Maisonobe and Merle's.

Let $(X,0)$ be a complex analytic germ defined by the vanishing of a
map of germs $F\:(\IC^n,0)\to (\IC^p,0)$ with $p\ge0$; if $p=0$, then
$F=\emptyset$ and $X=\IC^n$.  Let $f\:(\IC^n,0)\to (\IC,0)$ be the germ
of a complex analytic function.  Form the $p+1$ by $n$ matrix $D(F;f)$
by augmenting the Jacobian matrix $DF$ at the bottom with the gradient
$df$.  Call the submodule of the free module $\O_X^{p+1}$, generated by
the columns of $D(F;f)$, the {\it augmented Jacobian module\/} and
denote it by $JM(F;f)$.

More intrinsically, $JM(F;f)$ may be viewed as follows.  Identify
$(X,0)$ with the graph of $f|(X,0)$, which is a germ in $(\IC^{n+1},0)$.
This germ is defined by the vanishing of the map $G\:(\IC^{n+1},0)\to
(\IC^{p+1},0)$ whose components are $F$ and $f-z$, where $z$ is the last
coordinate function on $(\IC^{n+1},0)$.  Then $JM(F;f) = JM(G)_z$; that
is, they are the same submodule of $\O_X^{p+1}$.  Now, $JM(G)_z$ depends
only on the (abstract) normal module of the graph of $f|(X,0)$, not on
the choice of generators of this module (nor on the choice of
coordinates on $\IC^n$); see the beginning of Section~4.  Thus $JM(F;f)$
depends only on the restriction $f|(X,0)$; in other words, a second
function germ on $(\IC^n,0)$ with the same restriction as $f$ gives rise
to the same augmented Jacobian module, viewed as a submodule of the
normal module of the graph.  Moreover, given $F$ and $f$, this normal
module may be viewed as a submodule of the free module $\O_X^{p+1}$, and
the latter two modules are equal if $X$ is a complete intersection of
codimension $p$.

Given an analytic map germ $g\:(\IC^n,0)\to (\IC^l,0)$, let $JM(F;f)_g$
denote the submodule of $JM(F;f)$ generated by the columns of ``partial
derivatives'' with respect to the vector fields on $\IC^n$ tangent to
the fibers of $g$; in other words,
        $$JM(F;f)_g := JM(G)_{(g,z)},$$
 where $(g,z)\:(\IC^{n+1},0)\to (\IC^{l+1},0)$ has components $g$, $z$.
Call $JM(F;f)_g$ the {\it the relative augmented Jacobian module} with
respect to $g$.

 If $p=0$, or $F=\emptyset$, then write $JM(;f)$ and $JM(;f)_g$ for
$JM(F;f)$ and $JM(F;f)_g$.  Then $JM(;f)$ is simply the ideal on
$\IC^n$ generated by all the partial derivatives of $f$, and $JM(;f)_g$
is the subideal generated by the ``partials'' with respect to the
vector fields on $\IC^n$ tangent to the fibers of $g$.

Call a hyperplane in $\IC^n$ through $0$ a {\it limit tangent
hyperplane\/} of the fibers (or level hypersurfaces) of $f|X$ if it
is the limit of hyperplanes tangent to the fibers of $f|X$ at points
where $f|X$ is a submersion and that approach 0 along an analytic arc.
Now, let $(S,0)$ be a smooth subgerm of $(\IC^n,0)$ defined by the
vanishing of a map of germs $g\:(\IC^n,0)\to (\IC^l,0)$ with $l=n-m$
where $m:=\dim S$, and let $T_0S$ denote its tangent space at $0$.
Assume that $f|X$ is a submersion on the smooth locus of $X-S$.  The
following lemma describes the limit tangent hyperplanes of the fibers of
$f|X$ in general and those that contain $T_0S$ in particular in terms of
augmented Jacobian modules.

\lem1  Preserve the conditions above.
 \part1 A hyperplane $H$, defined by the vanishing of a linear
function $h\:\IC^n\to \IC$, is a limit tangent hyperplane of the fibers
of $f|X$ if and only if $JM(F;f)_h$ is not a reduction of $JM(F;f)$.
 \part2 No hyperplane containing $T_0S$ is a limit tangent hyperplane
of the fibers of $f|X$ if $JM(F;f)_g$ is a reduction of
$JM(F;f)$.
 \part3 Every limit tangent hyperplane of the fibers of $f|X$ contains
$T_0S$~--- in other words, the pair $(X-S,S)$ satisfies {\rm Thom's
Condition A}$_f$ at the origin~--- if and only if $JM(F;f)_g$ is
contained in $\sdep{JM(F;f)}$.
 \pf
 With $(F,f)$ in place of $F$, the proof is essentially the same as
that of Proposition~(4.1), because of the following observation: at a
nonsingular point $x$ of $X$, the Jacobian matrix of $(F,f)$ has
maximal rank because $x$ is not a critical point of $f|X$; moreover,
the row space of the matrix is the conormal module in $\IC^n$ to the
fiber through $x$ of $f|X$.

\bigbreak

Parameswaran \cite{P} generalized the L\^e--Ramanujam theorem
\cite{L-R, Thm.~2.1, p.~69} from a family of hypersurfaces to a family
of germs with isolated complete-intersection singularities (ICIS germs)
as follows.  To an ICIS germ $(X,0)$, he associated \cite{P, Def.~1,
p.~324} the sequence of numbers,
        $$\mu_*:=\mu_0,\mu_1,\dots,\mu_k,$$
 where $k$ is the embedding codimension and where $\mu_i$ is the
smallest Milnor number of any ICIS germ that serves as the total space
of a flat deformation of $(X,0)$ with a smooth parameter space of
dimension $i$.  Parameswaran noted \cite{P, Rmk., p.~324} that $\mu_k=0$
and that $\mu_i>0$ for $i<k$.  Given a chain (or nested sequence) of
deformations $(X_i,0)\to(Y_i,0)$ for $1\le i\le k$ such that each $Y_i$
is of dimension $i$, call the chain $\mu_*$-{\it minimal\/} if the
Milnor number of $(X_i,0)$ is equal to $\mu_i$.

Parameswaran proved \cite{P, Lem.~3, p.~325} that there exists a
$\mu_*$-minimal chain where each $Y_i$ is smooth.  He said \cite{P,
Def.~6, p.~331} that two ICISs have the same {\it topological type\/}
if each has a $\mu_*$-minimal chain such that the two chains are
embedded homeomorphic, and he proved that this notion does not depend
on the choice of chains.  Finally, Parameswaran proved \cite{P, Thm.~2,
p.~332} that, in a family of ICISs of dimension other than 2, the
topological types of the members are the same if their
$\mu_*$-sequences are the same; this is his generalization of the
L\^e--Ramanujam theorem.

Our next result stands to the L\^e--Saito theorem \cite{L-S, Thm.,
p.~793} as Parame\-swaran's result stands to the L\^e--Ramanujam
theorem; in fact, our result also asserts the converse to the
L\^e--Saito theorem, and generalizes it.  Consider a family of ICIS
germs, $(X,0)\to(Y,0)$ with section $\sigma (Y,0)\to (X,0)$ where $Y$ is
smooth, and let $k$ be the embedding codimension of $(X(0),0)$.
Parameswaran constructed a chain of deformations $(X_i,0)\to(Y_i,0)$ of
the family, with $1\le i\le k$, such that $Y_i/Y$ is smooth of relative
dimension $i$ and such that the fibers $(X_i(0),0)\to(Y_i(0),0)$ over
$0\in Y$ form a $\mu_*$-minimal chain (see the beginning of the proof of
\cite{P, Thm.~2, p.~332}).  Call such a chain a {\it Parameswaran chain}
if in addition, for each $y\in Y$ in a neighborhood of $0$, the fibers
$(X_i(y),0)\to(Y_i(y),0)$ over $y\in Y$ form a $\mu_*$-minimal chain.
Parameswaran also noted that the latter condition holds if the
$\mu_*$-sequence of $X(y)$ is constant in $y$ on a neighborhood of $0$
in $Y$.

Fix a Parameswaran chain.  Then $\mu_k=0$; so $(X_k,0)$ may be
identified with $(\IC^n,0)$ where $n:=\dim X_k$.  For convenience, set
$X_0:=X$ and $Y_0:=Y$.  For $0\le i< k$, let $f_i\:(X_k,0)\to(\IC^1,0)$
be a function that cuts $X_i$ out of $X_{i+1}$.  Finally, let $S$ denote
the singular locus of $X$.  Then, for $i<k$, the singular locus of $X_i$
is also $S$, and $f_i|X_{i+1}$ is a submersion off $S$.  Moreover, the
image of $Y$ under $\sigma$ lies in $S$, and $S$ is finite over $Y$.

\thm2 In the above setup, the pair $\bigl\leftp X_{i+1}-\sigma (Y),\sigma
(Y)\bigr)$ satisfies Thom's Condition {\rm A}$_{f_i}$ at the origin for
$0\le i<k$ if and only if the $\mu_*$-sequence of $(X(y),0)$ is constant
in $y$ on a neighborhood of $0$ in $Y$.
 \pf
 Let $h\:X_k\to Y$ be the projection.  For $0\le i<k$, set
        $$F_i:=(f_{k-1},\dots,f_{i+1})\and
                \M_i:=JM(F_i;f_i)_h\subset\O_{X_{i+1}}^{k-i},$$
 and let $J_i$ be the zeroth Fitting ideal of $\O_{X_{i+1}}^{k-i}/\M_i$.
Then, by the theorem of L\^e \cite{L, Thm.~3.7.1, p.~130} and Greuel
\cite{G, Kor.~,5.5, p.~263}, the colength of the induced ideal
$J_i\O_{X_{i+1}(y),0}$ is equal to the sum of the Milnor numbers of the
germs  $(X_{i+1}(y),0)$ and  $(X_i(y),0)$.  Hence the colength is
independent of $y\in Y$ near $0$ for all $i$
if and only if the Milnor numbers are so, since the Milnor numbers are
upper semicontinuous by \cite{Lo84, bot.~\p.126}.

Now, $X_{i+1}(y)$ has dimension $d+i+1$ where $d$ is the dimension of
$X(0)$, and $\M_i$ is generated by $d+k$ sections; hence, by virtue of
some theorems of Buchsbaum and Rim \cite{B-R, 2.4, 4.3, 4.5}, the
colength of the ideal $J_i\O_{X_{i+1}(y),0}$ is equal to the
Buchsbaum--Rim multiplicity, $e(i,y)$ say, of the image of $\M_i$ in
$\O_{X_{i+1}(y),0}^{k-i}$.  Therefore, since the chain is Parameswaran,
the $\mu_*$-sequence of $(X(y),0)$ is independent of $y$ near $0$ if and
only if, for each $i$, the multiplicity $e(i,y)$ is independent of $y$
near $0$.

The next part of the proof has some similarities with the beginning of
the proof of Theorem~(4.2).  Assume for the moment that $S=\sigma (Y)$.
Since $S$ is smooth, there exists a map germ $g\:(X_k,0)\to (\IC^l,0)$,
where $l:=\cod(S,X_k)$, such that $(S,0)=(g^{-1}0,0)$.  Moreover, by
\cite{Hir, Thm.~1, p.~242}, Thom's Condition A$_{f_i}$ is satisfied by
$(X_{i+1},S)$ at $s$ for all $s$ in some dense Zariski open subset
$U_i$ of $S$.  Hence, by Proposition~\Cs1), given $s\in U_i$,
        $$JM(F_i;f_i)_g\subset\sdep{JM(F_i;f_i)}\hbox{ at }s,\dno2.1$$
 and this relation holds at $s=0$ if and only if $(X_{i+1},S)$
satisfies A$_{f_i}$ at $0$.  Since the projection $S\to Y$ is an
isomorphism, the sum $JM(F_i;f_i)_g+\M_i$ is all of $JM(F_i;f_i)$.
Hence, by Part (i)$\Rightarrow$(ii) of Proposition~(3.2) and
Proposition~(3.1),
        $$JM(F_i;f_i)\?{}=\?{\M_i}\hbox{ at }s\dno2.2$$
 for $s\in U_i$, and this relation holds at $s=0$ if $(X_{i+1},S)$
satisfies A$_{f_i}$ at $0$.

Now, assume that the pair $\bigl(X_{i+1}-\sigma (Y),\sigma (Y)\bigr)$
satisfies A$_{f_i}$ at $0$ for each $i$; in particular, this assumption
means, by convention, that $\sigma (Y)$ is smooth, so $S=\sigma (Y)$.
Then, \Cs2.2) holds at $s=0$; so $\M_i$ is a reduction of $JM(F_i;f_i)$
over a neighborhood of $0$ in $Y$.  Hence the cosupport of $\M_i$ is
just $\sigma(Y)$.  Since $\M_i$ is generated by the right number of
sections, Proposition~(1.5)(3) says that $e(i,y)$ is independent of $y$
near 0, hence constant along $\sigma(Y)$.  Therefore, by the first
paragraph, the $\mu_*$-sequence of $(X(y),0)$ is independent of $y$
near $0$.

Conversely, assume that the $\mu_*$-sequence of $(X(y),0)$ is
independent of $y$ near $0$.  Then so is the multiplicity $e(i,y)$ for
each $i$ by the conclusion of the first paragraph.  Consider the
multiplicity of the image of $\M_i$ in $\O_{X_{i+1}(y)}^{k-i}$; it is
the sum of the multiplicities at each point of the fiber $X_{i+1}(y)$.
So it is at least $e(i,y)$, and the two are equal at $y=0$.  However,
the former is upper semicontinuous in $y$ by Proposition~(1.1).  Hence,
the two are equal for all $y$ near $0$.  Therefore, the cosupport of
$\M_i$ is equal to $\sigma(Y)$ over a Zariski open subset of $Y$, which
we may assume is all of $Y$.  Hence $S=\sigma (Y)$.  Hence $\M_i$ is a
reduction of $JM(F_i;f_i)$ over $h(U_i)$ because \Cs2.2) holds for
$s\in U_i$.  Therefore, Theorem~(1.8) implies that $\M_i$ is a
reduction of $JM(F_i;f_i)$.

Form the ideal $\rho(\M_i)$ on $X_{i+1}\x\IP^{k-i-1}$.  Form the
corresponding normalized blowup, its structure map,
  $$\pi\:NB_{\rho(\M_i)}(X_{i+1}\x\IP^{k-i-1})\to X_{i+1}\x\IP^{k-i-1},$$
 and its exceptional divisor $E$.  Since $\M_i$ is a reduction of
$JM(F_i;f_i)$, Proposition~(3.5)(1) yields the inclusion,
        $$\rho(JM(F_i;f_i)_g)\o\pi\subset\rho(\M_i)\o\pi.\dno2.3$$
 Moreover, this inclusion is strict along each component of $E$ that
projects onto $S$ because \Cs2.1) holds for $s\in U_i$.  Finally, to
complete the proof, it suffices, by Proposition~\Cs1), to prove that
the inclusion is strict along each remaining component $E_1$ of $E$.

Let $D$ be the exceptional divisor of the blowup itself of
$X_{i+1}\x\IP^{k-i-1}$ along $\rho(\M_i)$.  Then $E$ is the preimage of
$D$.  Since $\M_i$ has $d+k$ generators, $D$ lies in $\IP^{d+k-1}\x
S\x\IP^{k-i-1}$.  Now,
  $$\dim E_1 = \dim(X_{i+1}\x\IP^{k-i-1}) -1 = d+k-1+\dim Y,$$
 and $\dim S = \dim Y.$
 Also, $E_1$ does not project onto $S$.  Hence $E_1$ cannot project
onto a point of $\IP^{k-i-1}$.  Therefore, if $k=1$, then no such $E_1$
can exist, and the proof is complete in this case.

The proof proceeds by induction on $k$.  Let $L$ be a general
hyperplane in $\IP^{k-i-1}$.  Then the intersection
         $$\IP^{d+k-1}\x X_{i+1}\x L\bigcap
                B_{\rho(\M_i)}(X_{i+1}\x\IP^{k-i-1})$$
 is equal to the blowup $B$ of $X_{i+1}\x L$ along the ideal $\rho$ induced
by $\rho(\M_i)$.  Then $B$ contains the image $b$ of a general point of
$E_1$ because $L$ does.  Choose a map $\beta\:(\IC,0)\to(B,b)$ whose
image does not lie entirely in the exceptional divisor, and let
        $$\vf\:(\IC,0)\to(X_{i+1},0)\and
        \psi\:(\IC,0)\to(\Hom(\IC^{k-i},\IC),\lambda)$$
 be the maps arising from the composition of $\beta$ and the structure
map $B\to X_{i+1}\x\IP^{k-i-1}$.  If \Cs2.3) is not strict along $E_1$,
then
        $$\psi(\rho(JM(F_i;f_i)_g)\o\vf)=\psi(\M_i\o\vf).\dno2.4$$

Since $L$ is general, it is spanned by $k-i-1$ general points.  These
points correspond to $k-i-1$ general linear combinations
$g_{k-2},\dots,g_{i}$ of the $k-i$ functions $f_{k-1},\dots,f_i$, in
fact, to combinations of the functions cutting the $Y_i$ out of the
$Y_{i+1}$.  These functions define a chain of deformations $X_j'/Y_j'$
for $0\le j\le k-1$ where $X_{k-1}'$ and $Y_{k-1}'$ are $X_k$ and
$Y_k$.  For $j<k-1$, the singular locus of $X_j$ is $S$, and
$f_j|X_{j+1}$ is a submersion off $S$.  Moreover, by the upper
semicontinuity of Milnor numbers, the chain is Parameswaran, and the
$\mu_*$-sequence of $(X_i'(y),0)$ is independent of $y$ in a
neighborhood of $0$ in $Y$.  Thus the induction hypothesis applies.
Set $G_{i+1}:=(g_{k-1},\dots,g_{i+2})$.  Then \Cs2.4) becomes
        $$\psi(\rho(JM(G_{i+1};g_{i+1})_g)\o\vf)=
        \psi(\rho(JM(G_{i+1};g_{i+1})_h)\o\vf).$$
 However, this equation contradicts the induction hypothesis; indeed,
thanks to Lemma~(3.3) and Proposition~(5.1), the equation implies that
A$_{g_{i+1}}$ is not satisfied by the pair $(X'_{i+1}-S,S)$ at $0$.
The proof is now complete.
 \bigbreak

Brian\c con, Maisonobe and Merle found a relation between Whitney's
Condition A and Thom's Condition \Af, while working at the level of the
total space \cite{BMM, Thm.~4.2.1, p.~541}.  In essence, they proved
this.  Consider a pair $(X,Y)$ consisting of an analytic subspace $X$
of $\IC^n$, and a linear subspace $Y$ contained in $X$.  Consider a
function germ $f\:(\IC^n,0)\to (\IC,0)$.  Set
        $$Z:=f^{-1}(0)\cap X,$$
 and assume that $Z$ contains $Y$.  If both $X-Y$ and $Z-Y$ are smooth,
if both pairs, $(X-Y,Y)$ and $(Z-Y,Y)$, satisfy Whitney's Condition A
along $Y$, if $f$ is submersive on $X-Y$, and if, given the germ of any
linear retraction,
        $$r\:(\IC^n,0)\to(Y,0),$$
 the restriction $r|(X,Z,Y)$ is stratified locally topologically
trivial, then the pair $(X-Y,Y)$ satisfies Thom's Condition \Af\ along
$Y$.

Our next result, Theorem~(5.3) shows that, when $X$ is a complete
intersection and $Y$ is its singular locus, then the condition of
stratified triviality can be replaced by a numerical condition, which is
not only sufficient, but also necessary.  The numerical condition is
this: for every $r$ and for all $y\in Y$, the Buchsbaum--Rim
multiplicity,
        $$e(r,y):=e\bigl(JM(F;f)_{(r,y)}\bigr),$$
 is defined and constant in $y$, where $F\:(\IC^n,0)\to (\IC^p,0)$
defines $X$ as a complete intersection and where $JM(F;f)_{(r,y)}$
stands for the image of $JM(F;f)_r$ in $\O_{(r^{-1}(y)\cap X),0}^{p+1}$.

In Theorem~(5.3), it is unnecessary to assume that Whitney's Condition A
is satisfied.  Moreover, in Corollary~(5.4), we recover the original
theorem of Brian\c con, Maisonobe and Merle for ICIS germs in a refined
form: Condition A need be satisfied simply at 0, and the topological
trivializations of $r|(X,Y)$ and $r|(Z,Y)$ need not be compatible.

Furthermore, when Condition A is satisfied simply by $(Z-Y,Y)$ at 0,
then a much weaker condition will do.  It requires the constancy of
$e(r,y)$ only for a single $r$.  In other words, the condition depends
only on the individual fibers of $r|X$ and not on how they fit together
to form $X$.  See Theorem~(5.5).  In fact, Condition A is unnecessary
here too, as Massey and the first author proved in \cite{GaM, (5.8)}
after the present work was completed; see also \cite{Kl98} and
\cite{KT97, (1.7)}.

The Buchsbaum--Rim multiplicity $e(r,y)$ is defined if and only if the
ideal $JM(F;f)_{(r,y)}$ of $\O_{X,y}$ has finite colength; hence, since
$X$ is a complete intersection, if and only if the germs of the fibers
of the restriction $r|Z$ have isolated singularities.  If so, then the
germs of the fibers of the restriction $r|X$ have isolated singularities
too, and the following Milnor numbers are defined:
	$$\mu((r^{-1}(y)\cap X),0)\and\mu((r^{-1}(y)\cap Z),0).$$

These two Milnor numbers sum to $e(r,y)$ thanks to the theorem of L\^e
and Greuel and some theorems of Buchsbaum and Rim; see the first
paragraph of the proof of Theorem~\Cs2).  Since the two Milnor numbers
are upper semicontinuous by \cite{Lo84, bot.~\p.126}, they are
independent of $y$ if and only if $e(r,y)$ is so.  Hence we may
reformulate the three results below that involve the independence of
$e(r,y)$ by replacing this condition with the independence of the two
Milnor numbers.

\thm3
 In the setup of Brian\c con, Maxisonobe and Merle described above,
assume that $X$ is a complete intersection, and that both $X-Y$ and
$Z-Y$ are smooth.  Then the critical set $\Sigma(f)$ represents the same
germ as $Y$, and the pair $(X-Y,Y)$ satisfies \Af\ at $0$ if and only
if, for every linear retraction $r\:(\IC^n,0)\to(Y,0)$, the
Buchsbaum--Rim multiplicity $e(r,y)$ is defined and is independent of
$y$ for all $y\in Y$ near $0$.
 \pf
 To a certain degree, the proof is similar to those of Theorems~(4.2)
and (5.2).  Let $g\:(\IC^n,0)\to (\IC^l,0)$ be a map germ such that
$(Y,0)=(g^{-1}0,0)$ and $l:=\cod(Y,\IC^n)$.  Then, for all $y$ in some
dense Zariski open subset $U$ of $Y$, the pair $(X-Y,Y)$ satisfies \Af\
at $y$, and so
        $$JM(F;f)_g\subset\sdep{JM(F;f)}\hbox{ at }y\in U$$
  by Proposition~\Cs1)(3).  Hence, for every retraction $r$,
Proposition~(3.2) implies that ${JM(F;f)}_r$ is a reduction of
$JM(F;f)$ at $y\in U$.

Fix $r$, and assume that $e(r,y)$ is defined and independent of $y\in
Y$.  Consider the image of $JM(F;f)_r$ in $\O_{r^{-1}(y)\cap X}^{p+1}$.
It has finite colength for $y=0$, so for all $y$ in a neighborhood of
$0$, which we may assume is all of $Y$.  Since ${JM(F;f)}_r$ is
generated by the right number of sections, Proposition~(1.5)(3) yields
the constancy of its multiplicity.  This multiplicity is the sum of the
multiplicities at each point of the fiber $r^{-1}(y)\cap X$.  So it is
at least $e(r,y)$, and the two are equal at $y=0$.  Hence, the two are
equal for all $y$ near $0$.  Therefore, the cosupport of $JM(F;f)_r$ is
equal to $Y$.  Since this cosupport, $\Sigma(f)$, and $Y$ are always
nested, the three are equal.  In particular, $\Sigma(f)$ represents the
same germ as $Y$.

Therefore, Theorem~(1.8) implies that ${JM(F;f)}_r$ is a reduction of
$JM(F;f)$ everywhere.  Hence Lemma~\Cs1)(1) implies that no hyperplane
containing $\Ker r$ is a limit tangent hyperplane of the fibers of
$f|X$.  Now, given a hyperplane $H$ that does not contain $Y$, there
exists a retraction $r\:\IC^n\to Y$ such that $H$ contains $\Ker r$.
Therefore, $(X-Y,Y)$ satisfies Thom's Condition \Af\ at the origin.

Conversely, assume that $\Sigma(f)$ represents the same germ as $Y$, and
that the pair $(X-Y,Y)$ satisfies \Af\ at $0$.  Then $JM(F;f)_g$ is
contained in $\sdep{JM(F;f)}$ by Proposition~\Cs1)(3).  Hence, for every
retraction $r$, Proposition~(3.2) implies that ${JM(F;f)}_r$ is a
reduction of $JM(F;f)$ at $0$.  Since ${JM(F;f)}_r$ is generated by the
right number of sections, Proposition~(1.5)(3) yields the constancy of
$e(r,y)$ for all $y\in Y$ near $0$.  The proof is now complete.

\cor4 In the setup of Brian\c con, Maisonobe and Merle described
above, assume that $X$ is a complete intersection, that both $X-Y$ and
$Z-Y$ are smooth, and that both pairs $(X-Y,Y)$ and $(Z-Y,Y)$ satisfy
Whitney's Condition A at $0$.  Assume that, for every retraction $r$,
the restrictions $r|(X,Y)$ and $r|(Z,Y)$ are topologically trivial.
Then $(X-Y,Y)$ satisfies \Af\ at $0$.
 \pf We are about to prove that, after $X$ is replaced by a smaller
representative, the Milnor numbers, $\mu((r^{-1}(y)\cap X),0)$ and
$\mu((r^{-1}(y)\cap Z),0)$, are independent of $y$.  Hence, by the
discussion just before Theorem~(5.3), the Buchsbaum--Rim multiplicity
$e(r,y)$ is independent too.  Hence the assertion follows from
Theorem~\Cs3).

Since $(X-Y,Y)$ satisfies A at 0, it follows that, over a sufficiently
small neighborhood of $0$ in $Y$, the fibers of the restriction of $r|X$
are smooth except at points of $Y$.  Indeed, reasoning as in the last
paragraph of the proof of Theorem~(5.3), but using Lemma~(4.1)(3) in
place of Proposition~(3.2), we find that, on a sufficiently small
neighborhood, $JM(F)_r$ is a reduction of $JM(F)$.  Hence, the
cosupports of both these modules represent the same germ; otherwise,
there would be a map of germs $\vf\:(\IC,0)\to(X,0)$ whose image lies in
the former cosupport, but not in the latter, and then the pullbacks,
$JM(F)_r\o\vf$ and $JM(F)\o\vf$ would not be equal.  However, the
former cosupport is the singular locus $Y$ of $X$, and the latter
cosupport is the union of the singular loci of the fibers $r^{-1}(y)\cap
X$ of $r|X$.

Let $\Phi_y$ be the Milnor fiber of $r^{-1}(y)\cap X$, at its only
singular point $y$.  If $y$ is close enough to $0$, then there is a
short exact sequence of reduced integral homology groups,
	$$0\to\wt H_{n-p}(\Phi_y)\to\wt H_{n-p}(\Phi_0)
	       \to\wt H_{n-p}(r^{-1}(y)\cap X)\to 0;$$
 it is obtained from a versal deformation of $(X,0)$, see the top of
\p.121 in \cite{Lo84}.  Hence, if the Milnor number of $r^{-1}(0)\cap X$
is strictly greater than that of $r^{-1}(y)\cap X$, then the first
map cannot be surjective, and so $r^{-1}(y)\cap X$ is not contractible.

Replacing $X$ with a smaller representative, we may assume that
$r^{-1}(0)\cap X$ is contractible by \cite{Lo84, (2.4)}.  Then
$r^{-1}(y)\cap X$ is contractible too, since the restriction $r|(X,Y)$
is topologically trivial.  Therefore, the Milnor number of
$r^{-1}(y)\cap X$ is independent of $y$.  Similarly, the Milnor number
of $r^{-1}(y)\cap Z$ is independent of $y$ for $y$ near $0$, and the
proof is complete.

\thm5 In the setup of Brian\c con, Maisonobe and Merle described above,
above, assume that $X$ is a complete intersection, that both $X-Y$ and
$Z-Y$ are smooth, and that the $(Z-Y,Y)$ satisfies Whitney's Condition A
at $0$.  Then the following conditions are equivalent:\smallbreak
 \item i the critical set $\Sigma(f)$ represents the same germ as $Y$,
and the pair $(X-Y,Y)$ satisfies \Af\ at $0;$
 \item ii for every linear retraction $r\:(\IC^n,0)\to(Y,0)$, the
Buchsbaum--Rim multiplicity $e(r,y)$ is independent of $y$ for all
$y\in Y$ near $0;$
 \item ii$'$ for some linear retraction $r\:(\IC^n,0)\to(Y,0)$, the
Buchsbaum--Rim multiplicity $e(r,y)$ is independent of $y$ for all
$y\in Y$ near $0$.
 \pf
 Condition (i) implies (ii) by Theorem~\Cs3), and trivially (ii) implies
(ii$'$).  So assume (ii$'$), and let's prove (i).  Note that the proof
of Theorem~\Cs3) yields this: $\Sigma(f)$ represents the same germ as
$Y$, and no hyperplane containing $\Ker(r)$ is a limit tangent
hyperplane of the fibers of $f|X$.  Now, in $\IP^{n-1}$, the subspace of
hyperplanes containing $\Ker(r)$ has dimension $k-1$ where $k:=\dim Y$.
Since this subspace doesn't meet the space of limit tangent hyperplanes
of the fibers of $f|X$, the latter space must have dimension at most
$n-k-1$.

We now use an observation due to D. Massey and M. Green (pers.\ com.).
Form the relative conormal variety $C(X,f)$: by definition, it is the
closure in $\IC^n\x\IP^{n-1}$ of the locus of the pairs $(x,H)$ where
$x$ is a simple point of the level hypersurface surface $(f^{-1}fx)\cap
X$ and $H$ is a tangent hyperplane at $x$.  Intersect $C(X,f)$ with the
hypersurface $(f^{-1}0)\x\IP^{n-1}$.  Each component must have dimension at
least $n-1$ because $C(X,f)$ has dimension $n$.  Hence no component can
project onto a proper subset of $Y$; otherwise, the fiber of $C(X,f)$
over $0$ would have dimension at least $n-k$, but this fiber is simply
the space of limit tangent hyperplanes of the fibers of $f|X$, and
so it has dimension at most $n-k-1$ by the paragraph above.  Moreover,
if $(x,H)$ is a point of the intersection with $x\in Z-Y$, then $H$
must be tangent to $Z$ because $f$ is a submersion off $Y$ by
hypothesis.  Thus, each component of the intersection either surjects
onto $Y$ or lies in the conormal variety $C(Z)$; the latter is, by
definition, the closure of the locus of the pairs $(x,H)$ where $x$ is
a simple point of $Z$ and $H$ is a tangent hyperplane at $x$.

By hypothesis, $(Z,Y)$ satisfies Whitney's Condition \A\ at $0$; in
other words, the preimage of $Y$ in $C(Z)$ lies in $C(Y)$, the space of
hyperplanes containing $Y$.  Moreover, for all $y$ in a Zariski open
subset of $Y$, the pair $(X-Y,Y)$ satisfies \Af\ at $y$; in other words,
the fiber $C(X,f)(y)$ lies in $C(Y)$.  Hence, any irreducible subset of
$C(X,f)$ that projects onto $Y$ must lie in $C(Y)$.  Therefore, each
component of the intersection above lies in $C(Y)$.  So $C(X,f)(0)$
lies in $C(Y)$; in other words, $(X-Y,Y)$ satisfies $\A_f$.  The proof
is now complete.
 \medskip
 If $Y$ has dimension 1, then we have the following version of
Theorem~\Cs5).  It is a numerical criterion for Thom's Condition \Af,
which involves only a single retraction $r$ and not Whitney's Condition
A.

\cor6
  In the setup of Brian\c con, Maisonobe and Merle described above,
assume that $X$ is a complete intersection, that both $X-Y$ and
$Z-Y$ are smooth, and that $Y$ has dimension $1$.  Assume that, for each
hyperplane $H$ transverse to $Y$ at the origin, the Milnor numbers of
$H\cap Z$ and the Milnor number of a general hyperplane slice are
independent of $H$.  Then the following conditions are
equivalent:\smallbreak
 \item i the critical set $\Sigma(f)$ represents the same germ as $Y$,
and the pair $(X-Y,Y)$ satisfies \Af\ at $0;$
 \item ii for some linear retraction $r\:(\IC^n,0)\to(Y,0)$, the
Buchsbaum--Rim multiplicity $e(r,y)$ is independent of $y$ for all
$y\in Y$ near $0;$
 \item ii$'$ for every linear retraction $r\:(\IC^n,0)\to(Y,0)$, the
Buchsbaum--Rim multiplicity $e(r,y)$ is independent of $y$ for all
$y\in Y$ near $0$.
 \pf The hypothesis on the Milnor numbers is exactly what's needed to
conclude by Corollary~3.9 of \cite{G-3} that Whitney's Condition A
holds for the pair $(Z-Y,Y)$ at 0.  The result now follows from
Theorem~(5.5).

Here is the idea behind the proof of Corollary 3.9 of \cite{G-3}.
 Because $Y$ has dimension 1, a hyperplane $H:h=0$ that does not
contain $Y$ can intersect $Y$ only at  $0$.  Now, \cite{G-3}
makes a study of the hyperplane sections of $Z$ at $0$ by
hyperplanes transverse to $Y$.  On the basis the principle of
specialization of integral dependence, Theorem~(1.8) above, it is shown
in Proposition~2.7 and Theorem~3.3 of \cite{G-3} that the hypothesis on
the Milnor numbers implies that the submodule $JM(F,f)_h$ is a
reduction of $JM(F,f)$ in $\O_{ Z,0}^{p+1}$.  Hence Lemma~(4.1)(3)
implies that $H$ is not a limiting tangent hyperplane to $Z$ at 0.

\sctn The relative condition \Wf

In this section, given a map germ $f$ on $(X,0)$, we study the condition
\Wf.  It is a standard relative form of Whitney's Condition B, and
reduces to B, in the form of Verdier's Condition W \cite{Verdier,
Sect.~1}, when $f$ is constant.  In our first result, $X$ and $f$ are
arbitrary, but then we begin specializing as more hypotheses are needed.
In fact, we proceed to observe that $\Wf$\ is a rather strong condition
unless $f$ is a function germ, and from then on, we assume that $f$ is a
function.  Finally, in our last two result, we assume that $X$ is the
total space of a family of ICIS germs.

Condition \Wf\ generalizes Teissier's condition of `c-equisingularity'
(see \cite{LT88, top, p.~550}).  It strengthens Thom's Condition \Af,
and so is sometimes called the {\it strict Thom condition}.  Although
\Wf\ is defined using Euclidean distances, we prove in Proposition~\Cs1)
that \Wf\ is equivalent to a condition of integral dependence on a
modified Jacobian module, obtained by ``vertical'' differentiation.
Thus this module becomes the natural source for numerical invariants
that depend only on the members of a family $X/Y$, rather than on the
total space $X$.

The Thom--Mather second isotopy lemma readily implies that, if $f$ is a
nonconstant function and if \Wf\ is satisfied, then the pair $X,f$ is
topologically right trivial over $Y$.  Indeed, Thom, Mather, Teissier,
Verdier, and others introduced and developed methods of integrating
vector fields that yield this triviality.  Namely, we can lift to $X$ a
constant vector field tangent to $Y$ so that the lift is corrugated
(Fr.~{\it rugueux}), hence integrable, and is tangent to the fibers of
$f$ on $X/Y$, so that the integral gives a continuous flow on $X$.  If
we choose the field carefully, we can show that, after $X$ is replaced
by a neighborhood of 0, there is a homeomorphism $h\:X(0)\x Y\to X$ such
that $fh=(f|X(0))\x1_Y$, as required.  Similarly, it is possible to
generalize the statements and proofs of Corollary 3.6 and Theorem 3.8 of
\cite{G-2}.

Our main result, Theorem~\Cs4), characterizes \Wf, when $X$ is a family
of ICIS germs and $f$ is a nonconstant function, in three ways: (1) by
the constancy of the Buchsbaum--Rim multiplicity of a modified Jacobian
module, (2) by the constancy of two sequences of Milnor numbers, and (3)
by the fulfillment by two pairs, of the absolute Whitney conditions.
The necessity of (3) is trivial; its sufficiency is not new, but was
established by Brian\c con, Maisonobe and Merle in \cite{BMM,
Thm.~4.3.2, p.~543} in a more general setting using a different
approach.

We prove the theorem using Proposition~\Cs1) and Lemma~\Cs3).  The
latter expresses the Buchsbaum--Rim multiplicity in (1) as the weighted
sum of the Milnor numbers.  This lemma is proved using the polar
multiplicity formula and the relative polar transversality result of
Henry and Merle.  The latter is given a new proof in Lemma~\Cs2), and
this proof illustrates, for a second time, the usefulness of
Proposition~\Cs1) and of the methods of integral dependency.

A lemma similar to Lemma~\Cs3) is involved implicitly in the proof of
one of the main results, Theorem~1, in \cite{G-4}.  That proof does not
rely on the principle of specialization of integral dependence, our
Theorem~(1.8); indeed, the principle had not yet been established.
However, the principle yields a new proof of the implication
(iv)$\Rightarrow$(iii) of \cite{G-4, Thm.~2}, and this proof is in the
spirit of Teissier's original proof \cite{T-1} for the case of a
hypersurface.  In fact, the new proof is simply a special case of the
first part of the proof of the implication (iv)$\Rightarrow$(i) of our
Theorem~\Cs4); it is the case where $f$ is constant (so vanishes).  On
the other hand, (as the referee pointed out), it is possible to prove
this implication in the spirit of \cite{G-4}, using ordinary
multiplicities of polar varieties.

\medbreak

To begin the formal discussion, fix a pair $(X,Y)$ consisting of a
reduced equidimensional analytic subspace $X$ of $\IC^n$ and a linear
subspace $Y$ of $\IC^n$ contained in $X$.  Assume $l>p+q$.  Fix a map
germ $f\:(\IC^n,0)\to(\IC^q,0)$ whose restriction $f|(Y,0)$ is a
submersion onto a smooth closed analytic subgerm of $(\IC^q,0)$, and
assume that there is a smooth, dense, and open analytic subset $X_0$ of
$X$ such that $f|(X_0,0)$ is a submersion onto its image and has
equidimensional fibers.

 Recall from Definition~1.3.7 on \p.550 in \cite{LT88} (compare
\cite{H-M-S, pp.~228--9}) that $(X_0,Y)$ satisfies the {\it condition\/}
\Wf\ {\it at\/} $0$ if there exist a (Euclidean) neighborhood $U$ of $0$
in $X$ and a constant $C>0$ such that, for all $y$ in $U\cap Y$ and all
$x$ in $U\cap X_0$, we have
   $$\dist\bigl\leftp T_yY(f(y)), T_xX(f(x))\bigr)\leq C\,\dist(x,Y)$$
 where $T_yY(f(y))$ and $T_xX(f(x))$ are the tangent spaces to the
indicated fibers of the restrictions $f|Y$ and $f|X$.  This condition
depends only on the restrictions $F|X$ and $f|X$, and not on the
embeddings of $X$ into $\IC^n$ and of $f(X)$ into $\IC^q$.

Conditions like \Wf, which are defined by analytic inequalities, often
can be re-expressed algebraically in terms of integral dependence.  For
\Wf\ itself, this job was done by Navarro in a 1980 unpublished
manuscript according to Remarque 1.2(c) on \p.229 of \cite{H-M-S}.
Later the job was done in print by L\^e and Teissier.  In
Proposition~1.3.8 on \p.550 of \cite{LT88}, they translated \Wf\ into a
condition of integral dependence between ideals on the relative conormal
variety $C(X,f)$, whose definition was recalled in the proof of
Theorem~(5.5).  We recover their result below in Proposition~(6.1).

Proposition~(6.1) also gives another condition of integral dependence
equivalent to \Wf, and this is the condition of importance to us here.
It is the condition mentioned above, requiring that one modified
Jacobian module be dependent on another.  Before we can state and prove
the proposition formally, we must define these modules precisely.

Say that $(X,0)$ is defined by the vanishing of $F\:(\IC^n,0)\to
(\IC^p,0)$.  Generalizing the constructions in Section~5 , form the
corresponding {\it augmented Jacobian module} $JM(F;f)$: namely, first
form the $p+q$ by $n$ matrix $D(F;f)$ by augmenting the Jacobian of $F$
at the bottom with the Jacobian of $f$; then $JM(F;f)$ is the
$\O_X$-submodule of the free module $\O_X^{p+q}$, generated by the
columns of $D(F;f)$.

Say $\IC^n=\IC^l\x Y$, and form the corresponding projections,
	$$r\:\IC^n\to Y \and g\:\IC^n\to\IC^l.$$
 Form the corresponding {\it relative augmented Jacobian modules},
	$$JM(F;f)_r\and JM(F;f)_g;$$
 by definition, these are the submodules of $JM(F;f)$ generated by the
partial derivatives with respect to the first $l$ variables on $\IC^n$
and with respect to the remaining $n-l$ variables.  Finally, let $\mY$
be the ideal of $Y$ in $\IC^n$.

The abstract module $JM(F;f)$ is determined as a quotient of $O_X^n$,
but not as a submodule of $O_X^{p+q}$, by the (germ of the) embedding of
$X$ in $\IC^n$ and by the restriction $X\to S$, where $S$ is the image
$f(X)$ viewed as an abstract space.  Namely, $JM(F;f)$ is the unique
torsion free quotient that restricts to the normal sheaf on $X_0$.  So
the submodules $JM(F;f)_g$ and $JM(F;f)_r$ too are determined abstractly
by $X\to S$, given the splitting $\IC^n=\IC^l\x Y$.  Finally, it is
clear from its definition that the relative conormal variety $C(Y,f)$ is
determined by the embedding of $X$ in $\IC^n$ and by the restriction
$X\to S$ of $f$.

\prp1 In the setup above, the following conditions are equivalent:
 \smallskip
 \item i the pair $(X_0,Y)$ satisfies \Wf\ at $0;$
 \item ii the module $JM(F;f)_g$ is integrally dependent on $\mY
JM(F;f)_r$;
 \item ii$'$ the module $JM(F;f)_g$ is integrally dependent on
$\mY JM(F;f)$;
 \item iii along the preimage in $C(X,f)$ of $0$, the ideal of
   $C(Y,f)\cap C(X,f)$ is integrally
dependent on the ideal of the preimage of $Y$.
 \pf We'll prove that (ii$'$) is equivalent to each of the other
conditions.  First consider the notion of integral dependency involved
in (ii) and (ii$'$); it is defined abstractly in Section~1, but it can be
treated as discussed in Section~3, using the embedding of $JM(F;f)$ in
$O_X^{p+q}$.

That (ii) implies (ii$'$) is trivial.  Conversely, assume (ii$'$).  Then
$JM(F;f)_g$ is contained in the strict closure $\sdep{JM(F;f)}$.  Hence
$JM(F;f)$ is integrally dependent on $JM(F;f)_r$ by Prop.~(3.2) because
$JM(F;f)$ is the sum of $JM(F;f)_r$ and $JM(F;f)_g$.  So $\mY JM(F;f)$
is integrally dependent on $\mY JM(F;f)_r$.  Thus (ii$'$) implies
(ii).

To prove the equivalence of (i) and (ii$'$), let $e_1,\dots,e_n$ be a
vector space basis of $\IC^n$, and $f_1,\dots,f_{n-l}$ one of $Y$.  Then
the matrix products $D(F,f)\cdot e_i$ generate $JM(F;f)$, and the
products $D(F,f)_g\cdot f_j$ generate $JM(F;f)_g$.  Let $y_1,\dots,y_l$
be a set of coordinate functions on $\IC^l$.  Then the products
$y_kD(F,f)\cdot e_i$ generate $\mY JM(F;f)$.  So Proposition~1.11
on \p.306 of \cite{G-2} says that  (ii$'$) holds if and only if the
following condition holds: there exist a neighborhood $U'$ of $0$
in $X$ and a constant $C'>0$ such that, for any
$\psi\:U'\to\Hom(\IC^{p+q},\IC)$ and any $x\in U'$, we have
	$$\sup\nolimits_j|\psi(x)\cdot D(F,f)_g(x)\cdot f_j|
 \le C' \sup\nolimits_{i,k}|y_k(x)\psi(x)\cdot D(F,f)(x)\cdot e_i|.$$

The sup on the right is equal to
	$$\sup\nolimits_k|y_k(x)|\,
	\sup\nolimits_i|\psi(x)\cdot D(F,f)(x)\cdot e_i|.$$
 Adjusting the constant $C'$, we may replace the inequality above  by
 $$\|\psi(x)\cdot D(F,f)_g(x)\|\le C'\dist(x,Y)\|\psi(x)\cdot D(F,f)(x)\|.$$
 Set $u:=\psi(x)\cdot D(F,f)(x)$.  Then this inequality holds if and
only if, for every unit vector $v$ in $T_yY(f(y))$, the following
inequality holds:
	$$|(u,v)|\le C'\dist(x,Y)\|u\|.\dno1.1$$
 Here, we may replace $u$ by its complex conjugate.

If (i) holds, then the preceding inequality (6.1.1) holds with $U':=U$
and $C':=C$, at least for an $x$ in $U\cap X_0$, because, by definition,
	$$\dist(A,B):=\sup\nolimits_{\scriptstyle
	u\in B^{\perp}-\{0\}\atop\scriptstyle v\in A-\{0\}}
	{|(u,v)| \over \|u\|\|v\|}.$$
 By continuity, the inequality (6.1.1) also holds for an $x$ in $X-X_0$,
because the latter set is nowhere dense in $X$.  Thus (i) implies
(ii$'$).

Conversely, (ii$'$) implies (i).  Indeed, given any $x\in U\cap X_0$ and
$u\in B^{\perp}-\{0\}$ where $B:=T_xX(f(x))$, there is a
$\psi\:U'\to\Hom(\IC^{p+q},\IC)$ such that $\psi(x)$ is equal to the
conjugate of $u$; so we may take $U:=U'$ and $C:=C'$.

Finally, the equivalence of (iii) and (ii$'$) follows from the version
of Proposition~(3.4) for integral dependence given in Remark (10.8)(ii)
on \p.229 of \cite{K-T}.  In the latter, $\E$ is not necessarily locally
free, but we take
	$P:=\Projan(\R\E)$.
 (The proof is entirely different, and does not involve any form of the
valuative criterion.)  In the case at hand, take $\E$ to be $JM(F;f)$.
Then $P$ is just $C(X,f)$.

The ideal of the preimage in $P$ of $Y$ is just $\rho(\mY JM(F;f))$
because $\mY$ is the ideal of $Y$.  Furthermore, $\rho(JM(F;f)_g)$ is
the ideal of $C(Y,f)\cap C(X,f)$, because a hyperplane $\{w=0\}$ of
$\IC^n$ contains $Y$ if and only if the coefficients of the last $n-l$
coordinate functions of w are zero, and because the functions
corresponding to these coefficients are given on $C(X,f)$ by the columns
of $D(F,f)_g$.  The proof is now complete.
 \medbreak

Previous work on \Wf\ has involved the central fiber of the exceptional
divisor of the blowup of the relative conormal variety $C(X,f)$ along
the preimage of $Y$.  For example, Henry, Merle and Sabbah proved in
Cons\'equence~2 on \p.234 of \cite{H-M-S} that, if this fiber is of
minimal dimension and if \Wf\ holds generically on $Y$, then it also
holds at 0.  This conclusion also follows from Proposition~(6.1);
indeed, if the fiber is of minimal dimension and if (iii) holds
generically, then it also holds at 0 by B\"oger's celebrated criterion
of integral dependence of ideals (see \cite{K-T, (10.9)} and \cite{KT97,
(1.4)} for the generalization of this criterion to modules).
Conversely, assume that $f$ is a function or assume the more general
condition of the ``absence of blowup in codimension 0'' of Section~4 of
\cite{H-M-S}.  Then the fiber is of minimal dimension if \Wf\ holds
everywhere on $Y$.  This converse follows from Th\'eor\`eme~6.1 on
\p.262 and Proposition~3.3.1 on \p.239 of \cite{H-M-S}.  Recently, the
first author found a new proof using generic plane sections and methods
of integral dependence; the details will appear elsewhere.

Unless $f$ is a function (as it will be in our remaining three results),
\Wf\ is a very strong condition.  It implies that $f$ is analytically
right trivial already in this case: $X$ is $\IC^n$, the critical set
$\Sigma(f)$ of $f$ is reduced and is defined by the maximal minors of
the Jacobian matrix $D(f)$, and, for all $y\in Y$, the restriction of
the map germ $f(y):(\IC^l,0)\to(\IC^q,0)$ to its critical set is a
finite map onto its discriminant.  Indeed, say $f$ is nontrivial.  Then
both these latter sets have dimension $q-1$.  Moreover, as $y$ varies,
the union $\Sigma_Y(f)$ of these critical sets is equal to $\Sigma(f)$,
because $\Sigma_Y(f)$ is the cosupport of $JM(;f)_r$ and $\Sigma(f)$ is
the cosupport of $JM(;f)$; furthermore, the second module is integrally
dependent on the first by Proposition~(6.1).  (In fact, (6.1) is
stronger than necessary, and (5.1) will do after it is generalized from
a function to a map, a straightforward job; thus, already \Af\ implies
the analytic triviality of $f$.)

 Consider the map germ $F:=(f,r)$, with target $(\IC^q\x Y^k,0)$.  The
critical set of $F$ is just $\Sigma_Y(f)$, and we just proved that the
latter is equal to $\Sigma(f)$.  Since $F$ is finite on $\Sigma(f)$, its
discriminant, $\Delta(F)$, is a set of codimension 1. Now, $\Delta(F)$
projects onto the discriminant $\Delta(f)$, which is a proper subset of
$(\IC^q,0)$.  Hence $\Delta(F)$ is equal to $\Delta(f)\x Y$.  Now,
$\Sigma(F)$ is smooth of dimension $k+q-1$ on a dense Zariski open
subset $U$.  Shrinking $U$ if necessary, on $U$, the rank of $D(F)$ is
$k+q-1$, again because $F$ is finite on its critical set.  Hence, $\Ker
D(F)$ is transverse to $\Sigma(f)$ on $U$.  Therefore, because maximal
minors define $\Sigma(f)$ with reduced structure, $F$ is the unfolding
of a Morse function at points of $U$.  Hence, by Theorem ~1 on \p.726 of
\cite{BPW}, $F$ is analytically right trivial.

 \medbreak

Next we turn to the transversality result.  Let $P$ be a linear space
through $0$ in $\IC^n$ of codimension $i$ say, with $i\le\dim X$, and
let $\Pi$ be the {\it relative polar variety\/} of $f|X$ with $P$ as
pole.  By definition, $\Pi$ is the closure in $X$ of the locus of simple
points $x$ of the level hypersurface surface $X(f(x))$ such that there
exists a hyperplane that is tangent to $X(f(x))$ at $x$ and that
contains $P$.  In other words, $\Pi$ is the projection to $X$ of the
preimage $\nu^{-1}P^*$ where $\nu\:C(X,f)\to\IP^{n-1}$ is the projection
and $P^*$ is the set of hyperplanes containing $P$.  So, if $P$ is
general, then $\Pi$ has dimension $i-1+\dim S$.

Let $\pi\:\IC^n\to\IC^i$ be a linear map with kernel $P$.  Assume $f$ is
a nontrivial function, and let $\Sigma(f)$ be the critical locus of $f$
(which includes the singular locus of $X$).  Then $\Pi\cup \Sigma(f)$ is
cut out of $X$ by the maximal minors of the Jacobian matrix of the map
$\IC^n\to \IC^p\x\IC^q\x\IC^i$ with components $F$, $f$, and $\pi$.
Hence, if $P$ is general and if $\dim \Sigma(f)<i$, then $\Pi$ is
Cohen--Macaulay if $X$ is.

Remarkably, although $\Pi$ is defined using $P$, nevertheless the
two spaces are transverse at 0 if $P$ is general and $f$ is a
nonconstant function.  This is an important result.  Related results
were proved in the absolute case (the case where $f$ is constant) by
Teissier in \cite{T-1, 2.7--2.9}, in \cite{T75, Thm.~7, p.~623} and in
\cite{T77, Thm.~1, p.~269} for a hypersurface $X$, and by L\^e and
Teissier in \cite{LT81, (4.1.8), p.~569} for an arbitrary $X$.  Teissier
proved the relative result (where $f$ is a nonconstant function) for an
arbitrary $X$ in \cite{T81, pp.~40--41}, deriving it from his general
idealistic Bertini theorem.  This transversality result was also proved,
at about the same time, by Henry and Merle \cite{H&M82, Cor.~2, p.~195}.

The general relative polar transversality result is reproved next in a
new way, using the theory of the \Wf\ condition, especially
Proposition~(6.1) and the relative generic Whitney lemma.  The latter
was proved by Navarro, according to Henry, Merle, and Sabbah in Remarque
5.1.1 on \p.~255 of \cite{H-M-S}, and they generalized it (using the
normalized blowup of $C(X,f)$ along the preimage of $Y$ in the spirit of
Hironaka and of Teissier) in their Th\'eor\`eme~5.1 on the same page.

\lem2 (Relative polar transversality)\enspace In the setup above, the
relative polar variety $\Pi$ and its pole $P$ are transverse at 0 if $P$
is general and $f$ is a nonconstant function.
 \pf First note that the result is obvious if $i=\dim X$, as then
$\Pi=X$. Now, consider the ``Grassmann modification'' $\wt X$, which is
formed as follows.  Let $\IG$ be the Grassmann variety of all linear
spaces of codimension $i$ through $0$ in $\IC^n$, let $\wt\IC$ be the
tautological subbundle in $\IC^n\x\IG$, and let $\alpha\:\wt\IC\to\IC^n$
and $\beta\:\wt\IC\to\IG$ be the projections.  Set $\wt X:=\alpha^{-1}X$
and $\wt X_0:=\alpha^{-1}X_0$.  Then $\wt X_0$ is smooth since $X_0$ and
$\alpha$ are smooth.  Moreover, since $i<\dim X$, the 0-section
$0\x\IG$ of $\wt\IC$ lies in $\wt X$.

Set $\wt f:=f\circ\alpha$.  Then $\wt f|X_0$ is a submersion onto its
image and has equidimensional fibers because $\alpha$ is smooth with
equidimensional fibers.  Hence, since $\wt f$ is a function, by the
relative generic Whitney lemma, W$_{\wt f}$ is satisfied by $(\wt
X_0,0\x\IG)$ at $(0\x P)$ if $P$ lies in an appropriate dense Zariski
open subset of $\IG$.  Assume $P$ does so.

Proceeding via contradiction, assume that $P$ and $\Pi$ are not
transverse at 0.  Then there is a curve $\phi\:(\IC,0)\to(\Pi,0)$
tangent to $P$; moreover, in view of the definition of $\Pi$, we may
assume that, for $u\neq0$, we have $\phi(u)\in X_0$.  Since the
projection $C(X,f)\to X$ is proper, $\phi$ lifts to a map germ $\phi'$
from $(\IC,0)$ to $C(X,f)$ whose image lies in $\nu^{-1}P^*$.  So
$\nu\phi'(u)$ represents a hyperplane $H_u$ containing $P$ and, if
$u\neq0$, tangent to the fiber of $f|X$ through $\phi(u)$.  Let
$h\:\IC^n\to\IC$ be a linear functional whose kernel is equal to $H_0$.
Then $JM(F;f)_h$ is not a reduction of $JM(F;f)$ by Lemma~(5.1)(i); in
fact, the proof shows that $JM(F;f)_h\circ\phi$ is not equal to
$JM(F;f)\circ\phi$.  We'll now prove that they are equal since W$_{\wt f}$ is
satisfied; then we'll have a contradiction.

In $\IC^n$, choose an $i$-dimensional linear space $T$ through $0$ and
transverse to $P$.  Then the various $i$-codimensional spaces transverse
to $T$ form a Zariski neighborhood of $P\in\IG$, which we may identify
with $\Hom(\IC^i,\IC^{n-i})$.  In coordinates, a matrix $(a_{\mu,\nu})$
corresponds to the $i$-codimensional space with equations
	$t_\mu=\sum a_{\mu,\nu}z_\nu$
 where $t_1,\dots,t_i;z_1,\dots,z_{n-i}$ are coordinates on $\IC^n$
split as $T\x P$.  Over $\Hom(\IC^i,\IC^{n-i})$, the bundle $\wt\IC$ is
trivial, so equal to $P\x\Hom(\IC^i,\IC^{n-i})$.  In coordinates,
$\alpha\:\wt\IC\to\IC^n$ becomes
	$$\alpha(z_1,\dots,z_{n-i};a_{1,1},\dots,a_{k,n-i})
	 =(z_1,\dots,z_{n-i};
	\sum a_{1,\nu}z_\nu,\dots,\sum a_{k,\nu}z_\nu).$$
 Given a function germ $\gamma$ on $(\IC^n,0)$, note that we have
	$$\eqalignno{\pd{\gamma\circ\alpha}{a_{\mu,\nu}
		&=z_\nu\pd\gamma{t_\mu}\circ\alpha,&(6.2.1)\cr
	\pd{\gamma\circ\alpha}{z_\nu}&=\pd\gamma{z_\nu}\circ\alpha
 +\sum_\mu a_{\mu,\nu}\pd\gamma{t_\mu}\circ\alpha.&(6.2.2)\cr}}$$

Set $\wt F:=F\circ\alpha$.  Let $\delta$ denote the projection of
$P\x\Hom(\IC^i,\IC^{n-i})$ onto $P$, and let $\mG$ denote the ideal on
$\wt X$ generated by the $z_\nu$.  Since W$_{\wt f}$ is satisfied by $(\wt
X_0,0\x\IG)$ at $(0\x P)$, by Proposition~(6.1) the module
$JM(\wt F;\wt f)_\beta$ is integrally dependent on the product $\mG
JM(\wt F;\wt f)_\beta$.  Now, Equation~(6.2.1) implies that
	$$JM(\wt F;\wt f)_\delta=\mG (JM(F;f)_\tau\circ\alpha)$$
 where $\tau\:\IC^n\to P$ denotes the projection with kernel $T$.

Since the curve $\phi$ is tangent to $P$ at 0, there is a lift
$\wt\phi\:(\IC,0)\to(\wt X,0\x O)$.  Hence, by the preceding paragraph,
	$$(\mG(JM(F;f)_\tau\circ\alpha))\circ\wt\phi
	\subseteq (\mG JM(\wt F;\wt f)_\beta)\circ\wt\phi.$$
 The term on the left is equal to
$(\mG\circ\wt\phi)(JM(F;f)_\tau\circ\phi)$, and that on the right, to
$(\mG\circ\wt\phi)(JM(\wt F;\wt f)_\beta\circ\wt\phi)$.  So the
inclusion above is equivalent to this one:
	$$JM(F;f)_\tau\circ\phi\subseteq JM(\wt F;\wt f)_\beta\circ\wt\phi.$$
 Combined with Equation~(6.2.1), this inclusion implies the following
equation:
	$$JM(F;f)\circ\phi=JM(\wt F;\wt f)_\beta\circ\wt\phi.$$
 Now, Equation~(6.2.2) yields the inclusion,
	$$JM(\wt F;\wt f)_\beta\circ\wt\phi\subseteq
	  JM(F;f)_\zeta\circ\phi+{\bf m}(JM(F;f)_\tau\circ\phi),$$
 where $\zeta\:\IC^n\to T$ is the projection with kernel $P$ and where
$\bf m$ is the maximal ideal of $(\IC,0)$.  Hence Nakayama's lemma
yields the equation,
	$$JM(F;f)_\zeta\circ\phi=JM(\wt F;\wt f)_\beta\circ\wt\phi.$$
 The term on the right is equal to $JM(F;f)\circ\phi$ by the
preceding equation, and the term on the left lies in
$JM(F;f)_h\circ\phi$ since $P$ lies in $H$.  Therefore,
	$$JM(F;f)_h\circ\phi=JM(F;f)\circ\phi.$$
  Thus we have obtained the desired contradiction, and the proof is
complete.
 \medbreak

In our final two results, a key role is played by the level hypersurface,
	$$Z:=f^{-1}(0)\cap X.$$
 Assume  $Z\supset Y$, and assume $f|X$ is a {\it nonconstant\/}
function.  Given $y\in Y$, set
	$$X(y):=r^{-1}(y)\cap X\and Z(y):=r^{-1}(y)\cap Z.$$
 Finally, assume that $(X(y),0)$ and $(Z(y),0)$ are ICIS germs, and
that $(X,0)$ is given by the vanishing of $F\:(\IC^n,0)\to (\IC^p,0)$ as
a com\-plete intersection of codimension $p$.

Form the Buchsbaum--Rim multiplicity,
	$$em(y):=e(\mY\, JM(F;f)_r(y)),$$
  where $JM(F;f)_r(y)$ stands for the image of $JM(F;f)_r$ in
$\O_{X(y),0}^{p+1}$.  Since $\mY$ induces the maximal ideal ${\bf
m}_{X(y)}$ and $JM(F;f)_r(y)$ is equal to the augmented Jacobian module
$JM(F|X(y);f|X(y))$ of the restricted functions, we have
	$$em(y)=e({\bf m}_{X(y)}\,JM(F|X(y);f|X(y))),$$
 which is an invariant of the fiber over $y\in Y$.

\lem3 In the setup above, fix $y$, and for $i=0,\dots,l-p$, let
$\mu_i(X(y),0)$ and $\mu_i(Z(y),0)$ denote the Milnor numbers of the
sections by a general linear space $P_i$ of codimension $i$ in $\IC^l$.
Then
   $$em(y)=\Sum_{i=0}^{l-p}{l-1\choose i}\bigl\leftp\mu_i(X(y),0)
		+\mu_i(Z(y),0)\bigr).$$
 \pf
 Let $\Pi^i_y$ denote the $i$-dimensional relative polar subscheme of
$f|X(y)$ with $P_i$ as pole.  It follows from the polar multiplicity
formula \cite{K-T, Thm.~(9.8)(i)} (compare with \cite{H-M, 4.2.7} and
\cite{G-4, \S3}) that
	$$em(y)= \Sum_{i=0}^{l-p}{l-1\choose i}m(\Pi^i_y,0)$$
 where $m(\Pi^i_y,0)$ is the ordinary multiplicity at $0$ of $\Pi^i_y$.
In particular, $m(\Pi^{l-p}_y,0)$ is simply the multiplicity of $X(y)$
at $0;$ so it is equal to $\mu_{l-p}(X(y),0)+1$.  By convention,
$\mu_{l-p}(X_y,0)$ and $\mu_{l-p-1}(Z_y,0)$ and $\mu_{l-p-1}(Z_y,0)$ are
the ordinary multiplicities at $0$ diminished by $1$, and
$\mu_{l-p}(Z_y,0)=1$.

In general, since $\Pi^i_y$ is Cohen--Macaulay,
	$$m(\Pi^i_y,0) = \dim\O_{\Pi^i_y,0}/\I(L)\dno2.1$$
 where $\I(L)$ is the ideal of any linear space $L$ of codimension $i$ in
$\IC^l$ that is transverse to $\Pi^i_y$.  By the preceding lemma, we may
take $P_i$ for $L$.  Then, for $i<l-p$, the right side is equal to
$\mu_i(X(y),0)+ \mu_i(Z(y),0)$ by the theorem of L\^e and Greuel.  The
asserted formula follows immediately, and the proof is complete.

 \thm4 In the setup of Lemma~(6.3), let $\Sigma(f)$ denote the the
critical set of $f$, and $\Sigma_Y(f)$ the union of the critical sets
of the restrictions $f|X(y)$.  Then the following four conditions are
equivalent:\smallbreak
 \item i the germs of $\Sigma(f)$ and $Y$ are equal, and the pair
$(X-Y,Y)$ satisfies \Wf\ at $0;$
 \item ii the germs of $\Sigma_Y(f)$ and $Y$ are equal, and both
pairs $(X-Y,Y)$ and $(Z-Y,Y)$ satisfy the absolute Whitney conditions at
$0;$
 \item iii the Milnor numbers of the sections, $\mu_i(X(y),0)$ and
$\mu_i(Z(y),0)$, are constant in $y\in Y$ near $0;$
 \item iv the multiplicity $em(y)$ is constant in $y\in Y$ near $0$.
 \pf First of all, (i) implies (ii).  Indeed, $\Sigma_Y(f)$ and
$\Sigma(f)$ represent the same germ by the argument given in the third
paragraph before Lemma~(6.2).  Furthermore, $T_xX(f(x))\subset T_xX$
and, if $x\in Z$, then $T_xX(f(x))=T_xZ$.  Hence, the analytic
inequalities required by (ii) are automatically satisfied when (i)
holds.  Second, (ii) implies (iii); indeed, this implication is
virtually the assertion of Th\'eor\`eme~(10.1) on \p.223 of \cite{N80}.
Third, (iii) implies (iv) by Lemma~(6.3).

Lastly, assume (iv).  Then the germs of $\Sigma(f)$ and $Y$ are equal by
the upper semicontinuity argument in the proof of Theorem~(5.3), this
time applied to ${\bf m}_Y JM(F;f)_r$.  Now, since $f$ is a function, by
the relative generic Whitney lemma, \Wf\ holds at a general point of
$Y$.  Hence, generically $JM(F;f)_g$ is integrally dependent on $\mY
JM(F;f)_r$ by Proposition~(6.1).  Therefore, since (iv) holds, this
dependency holds at 0 by the principle of specialization of integral
dependence, Theorem~(1.8).  So Proposition~(6.1) implies that
$(X-\Sigma(f),Y)$ satisfies \Wf\ at $0$.  Thus (i) holds and the proof
is complete.

\sct References

\references

BMM
J. Brian\c con{,} P. Maisonobe and M. Merle,
 {\it Localisation de syst\`emes
diff\'erentiels, stratifications de Whitney et condition de Thom,}
 \invent 117 1994 531--50

BPW
 J. W. Bruce, A. A. Du Plessis, and L. C. Wilson,
 {\it Discriminants and liftable vector fields,\/}
 \jag 3 1994 725--53

B-R
 D. A. Buchsbaum and D. S. Rim,
  {\it A generalized Koszul complex. II. Depth and multiplicity,}
 \tams 111 1963 197--224

Fulton
  W. Fulton,
 ``Intersection Theory,''
 Ergebnisse der Mathematik und ihrer Grenzgebiete, 3. Folge
 $\cdot$ Band 2, Springer--Verlag, Berlin, 1984

G-1
 T. Gaffney,
 {\it Aureoles and integral closure of modules,}
 in ``Stratifications, Singularities and Differential ~II. Travaux en
Cours, {\b 55},'' Herman, Paris, (1997), 55--62.

G-2
 T. Gaffney,
 {\it Integral closure of modules and Whitney equisingularity,}
 \invent 107 1992 301--22

G-3
 T. Gaffney,
 {\it Equisingularity of plane sections, $t_1$ condition, and the
integral closure of modules,} in ``Real and Complex Singularities"
Proceedings of the Third International Workshop on Real and Complex
Singularities at Sao Carlos, Brasil 1994, W. L. Marar (ed.) Pitman
Research Notes in Mathematics 333 (1995) 95-111

G-4
 T. Gaffney,
 {\it Multiplicities and equisingularity of ICIS germs,}
\invent 123 1996 209--20

GaM
 T. Gaffney and D. Massey,
 {\it Trends in equisingularity theory,}
 to appear in the proceedings of a 1996 symposium in honor of CTC Wall,
 singularities volume, W. Bruce and D. Mond (eds.), Cambr. U. Press, to
appear.

Gr
 M. D. Green,
 {\it Dissertation, Northeastern University,} 1997

G-M
 M. D. Green and D. B. Massey,
 {\it Vanishing cycles and Thom's $a_f$ conditions,} Preprint 1996

G
 G. M. Greuel,
 ``Der Gauss--Manin Zusammenhang isolierter Singularit\"aten von
vollst\"and\-igen Durchschnitten," Dissertation, G\"ottingen (1973),
\ma 214 1975 235--66

H&M82
 J.P.G. Henry and M. Merle,
 {\it Limites d'espaces tangents et transversalit\'e de vari\'et\'es
polaires,}
 in ``Proc. La R\'abida, 1981.'' J. M. Aroca, R. Buchweitz, M.
Giusti and M.  Merle (eds.) \splm 961 1982 189--99

H&M83
 J.P.G. Henry and M. Merle,
 {\it Limites de normales, conditions de Whitney et \'eclatement d'Hironaka,}
 ``Singularities,''
 Proc. Symposia pure math. vol. {\bf 40}, part 1,
 Amer. Math.  Soc. (1983), 575--84

H-M
 J.P.G. Henry and M. Merle,
 {\it Conormal Space and Jacobian module.  A short dictionary,}
 in ``Proceedings of the Lille Congress of Singularities," J.-P. Brasselet
(ed.), London Math. Soc. Lecture Notes {\bf 201} (1994), 147--74

H-M-S
 J.P.G. Henry, M. Merle,  and C. Sabbah,
 {\it Sur la condition de Thom stricte pour un morphisme analytique
complexe,}
 \asens 17 1984 227--68

Hir
 H. Hironaka,
  {\it Stratification and flatness,} in ``Real and complex
singularities, Nordic Summer School, Oslo, 1976,'' Sijthoff \&
Noordhoff, (1977), 199--265

R-K
 D. Kirby and D. Rees,
{\it Multiplicities in graded rings I: The general theory,}
 in ``Commutative algebra: syzygies, multiplicities, and birational
algebra'' W. J. Heinzer, C. L. Huneke, J. D. Sally (eds.), \conm 159 1994
209--67

Kl98
 S. L. Kleiman
 {\it Equisingularity, multiplicity, and dependence,}
 to appear in the proceedings of a conference in honor of M. Fiorentini,
Marcel Dekker, 1998

K-T
 S. Kleiman and A. Thorup,
  {\it A geometric theory of the Buchsbaum--Rim multiplicity,}
 \ja 167 1994 168--231

KT97
 S. Kleiman and A. Thorup,
  {\it Conormal geometry of maximal minors,}
 alg-geom/970818

L
 D. T. L\^e,
  {\it Calculation of Milnor number of isolated singularity of complete
intersection,}
 \faa 8 1974 127--31

L-R
 D. T. L\^e and C. P. Ramanujam,
 {\it The invariance of Milnor's number implies the invariance of the
topological type,}
 \ajm 98 1976 67--78

L-S
 D. T. L\^e and K. Saito,
 {\it La constance du nombre de Milnor donne des bonnes stratifications,}
 \crasp 277 1973 793--95

LT81
 D. T. L\^e and B. Teissier,
 {\it Vari\'et\'es polaires locales et classes de Chern des
vari\'et\'es singuli\`eres}
 \ann 114 1981 457--91

LT88
 D. T. L\^e and B. Teissier,
 {\it Limites de'espaces tangent en g\'eom\'etrie analytique,}
 \cmh 63 1988 540--78

LJT
 M. Lejeune-Jalabert and B. Teissier,
 {\it Cl\^oture integrale des ideaux et equisingularit\'e, chapitre 1}  Publ.
 Inst. Fourier (1974)

Lip
 J. Lipman,
  {\it Equimultiplicity, reduction and blowing-up,}
  in ``Commutative algebra: analytic methods,''
 \dlnpam 68 1982 111--48

Lo84
 E.J.N.  Looijenga,
 Isolated singular points on complete intersections,
 London Mathematical Society lecture note series {\bf 77},
Cambridge University Press, 1984.

N80
 V. Navarro,
 {\it Conditions de Whitney et sections planes,}
 \invent 61 1980 199--226

P
 A. J. Parameswaran,
 {\it Topological equisingularity for isolated complete intersection
singularities,}
 \comp 80 1991 323--36

R87
 D. Rees,
 ``Reduction of modules,"
 \mpcps 101 1987 431--49

R89
 D. Rees,
 {\it Gaffney's problem,}
  manuscript dated Feb.22 1989

T-1
 B. Teissier,
 {\it Cycles \'evanescents, sections planes et conditions de Whitney,}
 in ``Singularit\'es \`a Carg\`ese," Ast\'erisque {\bf 7--8} (1973),
285--362

T75
B. Teissier,
{\it Introduction to equisingularity problems}
 Proc. Symposia pure math.  vol. {\bf 29}
Amer. Math.  Soc. (1975), 575--84

T76
 B. Teissier,
 {\it The hunting of invariants in the geometry of the discriminant,}
 in ``Real and complex singularities, Oslo
1976,'' P. Holm (ed.), Sijthoff \& Noordhoff (1977), 565--678

T77
B. Teissier,
{\it Vari\'et\'es polaires. I}
 \invent 40 1977 267--92

T1.5
 B. Teissier,
 {\it R\'esolution simultan\'ee et cycles \'evanescents,}
 in ``S\'em. sur les singularit\'es des surfaces.''  Proc.  1976--77.  M.
Demazure, H. Pinkham and B. Teissier (eds.) \splm 777 1980 82--146

T81
 B. Teissier,
 {\it Vari\'et\'es polaires locales: quelques r\'esultats},
 in ``Journ\'ees complexes,'' Institut Elie Cartan, Nancy, March 1981

T-2
 B. Teissier,
 {\it Multiplicit\'es polaires, sections planes, et conditions de
Whitney,}
 in ``Proc. La R\'abida, 1981.'' J. M. Aroca, R. Buchweitz, M. Giusti and
M.  Merle (eds.) \splm 961 1982 314--491

Trot
 D. Trotman,
  {\it On the canonical Whitney stratification of algebraic hypersurfaces,}
 Sem. sur la g\'eom\'etrie alg\'ebrique realle, dirig\'e par J.-J.
Risler, Publ. math. de l'universit\'e Paris VII, Tome I (1986), 123--52

Verdier
 J.-L. Verdier,
 {\it Stratifications de Whitney et th\'eor\`eme de Bertini--Sard,}
\invent 36  1976 295--312

\endreferences

\bigskip
    \parindent=\keyindent
 \eightpoint\smc
 Dept.~of Mathematics, Northeastern University, Boston, MA 02115, USA\par
 {\it E-mail address\/}: {\tt Gaff@NeU.edu}
 \smallskip
 Dept.~of Mathematics, 2--278 MIT, Cambridge, MA 02139, USA\par
 {\it E-mail address\/}: {\tt Kleiman@math.MIT.edu}

\bye